\documentclass[aop,preprint]{imsart}

\usepackage{amsmath, amsfonts, amsthm, verbatim, amssymb}

\newcommand{\dd}{\mathrm{d}}
\newcommand{\R}{\mathbb{R}}
\newcommand{\Z}{\mathbb{Z}}

\newcommand{\N}{\mathbb{N}}
\newcommand{\E}{\mathbb{E}}

\theoremstyle{plain}
\newtheorem{theorem}{Theorem}[section]
\newtheorem{corollary}[theorem]{Corollary}
\newtheorem{lemma}[theorem]{Lemma}
\newtheorem{proposition}[theorem]{Proposition}
\newtheorem{theoremA}{Theorem}

\theoremstyle{remark}
\newtheorem{remark}[theorem]{Remark}

\theoremstyle{definition}
\newtheorem{definition}[theorem]{Definition}
\newtheorem{example}[theorem]{Example}
\newtheorem{assumption}{Assumption}

\usepackage{latexsym}
\usepackage{amssymb}
\usepackage{amsmath}
\usepackage{amsthm}
\usepackage{mathrsfs}
\usepackage[numbers]{natbib}
\usepackage{url}
\usepackage{longtable}
\usepackage{color}
\usepackage{dsfont}
\usepackage[T1]{fontenc}
\usepackage[utf8]{inputenc}
\usepackage{lmodern}
\usepackage[colorlinks,citecolor=blue,urlcolor=blue]{hyperref}

\usepackage{graphicx}

\usepackage{amsopn}

\startlocaldefs

\newcommand{\bfi}{\begin{fig}}
	\newcommand{\efi}{\end{fig}}
\newcommand{\btab}{\begin{tab}}
	\newcommand{\etab}{\end{tab}}
\newcommand{\barr}{\begin{array}}
	\newcommand{\earr}{\end{array}}
\newcommand{\beq}{\begin{equation}}
\newcommand{\eeq}{\end{equation}}
\newcommand{\bdis}{\begin{displaymath}}
\newcommand{\edis}{\end{displaymath}\noindent}

\newcommand{\bbe}{\mathbb{E}}

\newcommand{\bbp}{\mathbb{P}}

\newcommand{\bone}{\mathds 1}

\newcommand{\limw}{\stackrel{\rm w}{\longrightarrow}}

\newcommand{\limq}{\stackrel{L^q}{\Longrightarrow}}

\newcommand{\lec}{\lesssim}

\newcommand{\calu}{{\cal U}}
\newcommand{\calf}{{\cal F}}

\newcommand{\call}{{\cal L}}

\newcommand{\calg}{{\cal G}}
\newcommand{\calq}{{\cal Q}}

\newcommand{\calo}{{\cal O}}

\newcommand{\al}{{\alpha}}
\newcommand{\la}{{\lambda}}

\newcommand{\eps}{{\varepsilon}}
\newcommand{\ga}{{\gamma}}
\newcommand{\Ga}{{\Gamma}}

\newcommand{\si}{{\sigma}}

\newcommand{\Om}{{\Omega}}

\newcommand{\var}{{\mathrm{Var}}}
\newcommand{\cov}{{\mathrm{Cov}}}

\newcommand{\ov}{\overline}

\newcommand{\wt}{\widetilde}

\newcommand{\ee}{\mathrm{e}}

\definecolor{orange}{rgb}{1,0.5,0}

\newcommand{\lb}{\left[}
\newcommand{\rb}{\right]}
\newcommand{\lp}{\left(}
\newcommand{\rp}{\right)}
\newcommand{\lv}{\left|}
\newcommand{\rv}{\right|}
\newcommand{\as}{\quad\text{as}\quad n\to\infty}
\newcommand{\sumt}{\sum_{i=1}^{[t/\Delta_n]}}

\newcommand{\Del}{\Delta_n}

\newcommand{\supt}{\sup_{t\in[0,T]}}

\DeclareMathOperator*{\dist}{dist}

\newcommand{\bthm}{\begin{theorem}}
	\newcommand{\ethm}{\end{theorem}}

\newcommand{\bthmA}{\begin{theoremA}}
	\newcommand{\ethmA}{\end{theoremA}}

\newcommand{\bcor}{\begin{corollary}}
	\newcommand{\ecor}{\end{corollary}}

\newcommand{\blem}{\begin{lemma}}
	\newcommand{\elem}{\end{lemma}}

\newcommand{\bprop}{\begin{proposition}}
	\newcommand{\eprop}{\end{proposition}}

\newcommand{\bdf}{\begin{definition}}
	\newcommand{\edf}{\end{definition}}

\newcommand{\bex}{\begin{example}}
	\newcommand{\eex}{\end{example}}

\newcommand{\brem}{\begin{remark}}
	\newcommand{\erem}{\end{remark}}

\newcommand{\bass}{\begin{assumption}}
	\newcommand{\eass}{\end{assumption}}

\newcommand{\bpr}{\begin{proof}}
	\newcommand{\epr}{\end{proof}}

\newcommand{\benu}{\begin{longlist}}
	\newcommand{\eenu}{\end{longlist}}

\newcommand{\bit}{\begin{itemize}}
	\newcommand{\eit}{\end{itemize}}

\newcommand{\bff}{\textbf}

\numberwithin{equation}{section}

\allowdisplaybreaks[4]
\endlocaldefs

\begin{document}
	
	\begin{frontmatter}
\title{Power variations in fractional Sobolev spaces for a class of parabolic stochastic PDEs}
		\runtitle{Power variations for SPDEs}
		
		\begin{aug}
			\author[A]{\fnms{Carsten} \snm{Chong}\ead[label=e1,mark]{carsten.chong@epfl.ch}}
			\and
			\author[A]{\fnms{Robert C.} \snm{Dalang}\ead[label=e2,mark]{robert.dalang@epfl.ch}}
		
			\address[A]{Institut de mathématiques, École Polytechnique Fédérale de Lausanne\\
				\printead{e1,e2}}
		\end{aug}
		
	\begin{abstract}
		We consider a class of parabolic stochastic PDEs on bounded domains $D\subseteq\mathbb{R}^d$  that includes the stochastic heat equation, but with a fractional power $\gamma$ of the Laplacian. Viewing the solution as a process with values in a scale of fractional Sobolev spaces $H_r$, with $r < \gamma - d/2$, we study  its power variations in $H_r$ along regular partitions of the time-axis.  As the mesh size tends to zero, we find a phase transition at $r=-d/2$:  the solutions have a  nontrivial quadratic variation when $r<-d/2$ and a nontrivial $p$th order variation for $p= 2\gamma/(\gamma-d/2-r)>2$ when $r>-d/2$.   More generally, suitably normalized power variations of any order satisfy  a genuine law of large numbers in the first case and a degenerate limit theorem in the second case. When $r<-d/2$, the quadratic variation is given explicitly via an expression that involves the spectral zeta function, which reduces to the Riemann zeta function when $d=1$ and $D$ is an interval.
		%
	\end{abstract}
		
		
		
	\end{frontmatter}

\section{Introduction}

 Let $D$ be a bounded open subset of $\R^d$ (satisfying certain regularity conditions) and consider the following parabolic stochastic PDE on $[0,\infty)\times D$ with zero Dirichlet boundary conditions:
\beq\label{SHE}
\left\{\begin{array}{ll} \frac{\partial u}{\partial t}(t,x) = -(-\Delta)^\ga u(t,x)+\sigma(u(t,x)) \dot W(t, x), & (t,x) \in [0,\infty)\times D,\\ u(t,x)=0, &  (t,x)\in[0,\infty)\times\partial D, \\
	u(0,x)=0, & x\in D. \end{array}\right.\eeq
Here, $\dot W$ is a Gaussian space-time white noise on $[0,\infty)\times D$, $\si\colon \R\to\R$ is a Lipschitz function, and $(-\Delta)^\ga$ is a spectral power of $-\Delta$ (see Section~\ref{sec2} for details). The purpose of this article is study the regularity of $t\mapsto u(t,\cdot)$ as a stochastic process taking values in a scale of Sobolev spaces $H_r=H_r(D)$ indexed by $r\in\R$ (to be defined in Section~\ref{sec2}). More precisely,  we are interested in the asymptotic behavior as $n\to\infty$ of the \emph{(normalized) power variation of order $p$} given by
\beq\label{Vnp} V^{n,r}_p(u,t):= \Delta_n\sum_{i=1}^{[t/\Delta_n]} \biggl( \frac{\|u(i\Delta_n,\cdot)-u((i-1)\Delta_n,\cdot)\|_{H_r}}{\tau_n(r)} \biggr)^p,\quad t\in[0,\infty),\quad n\geq1. \eeq
In this paper, $\Delta_n$ is a strictly positive sequence decreasing to $0$ (e.g., $\Delta_n=\frac1n$), $p>0$ is a fixed but arbitrary power, and $\tau_n(r)$ is a normalizing factor depending on $r$  and $n$ (and $\ga$), chosen if possible in such a way that $V^{n,r}_p(u,t)$ converges to a limit $V^r_p(u,t)$, say, uniformly on compact sets in probability.

In order to describe the flavor of our results, let us specialize to the case where $d=1$, $D=(0,\pi)$, $\ga=1$, and $\si\equiv1$ (i.e., to the stochastic heat equation on an interval with additive noise) in this introductory part. As Proposition~\ref{uinHr} below shows, in order that $u(t,\cdot)\in H_r$ for $t>0$, the smoothness parameter $r$ must be taken in the range $(-\infty,\frac12)$.
\bthmA\label{thmA} Let $u$ be the solution of \eqref{SHE} with $d=1$, $D=(0,\pi)$, $\ga=1$, and $\si\equiv1$. Assume $r<\frac12$ and define
\beq\label{taunr0} \tau_n(r):=\begin{cases} \Del^{\frac12} &\text{if } r<-\frac12,\\ (\Del \lvert\log\Del\rvert)^{\frac12} & \text{if } r=-\frac12,\\ \Del^{\frac14-\frac r2} &\text{if } -\frac12<r<\frac12. \end{cases} \eeq
If $p\geq1$ is an integer, then
\beq\label{conv} \lim_{n\to\infty} \E\Biggl[ \sup_{t\in[0,T]} \Bigl|V^{n,r}_{2p}(u,t)-K(r,p)t\Bigr|^q\Biggr] =0 \eeq
for every $q,T\in(0,\infty)$, where the constant $K(r,p)$ is given by
\beq\label{const} K(r,p):=\begin{cases} 2^p B_p\Bigl( \frac{(1-1)!}{2}\zeta(-2\cdot  1 \cdot r),\dots,\frac{(p-1)!}{2}\zeta(-2\cdot  p \cdot r) \Bigr)&\text{if } r<-\frac12,\\ 2^{-p}& \text{if } r=-\frac12,\\ \biggl( \frac{\Ga(r+\frac12)}{2(\frac12-r)} \biggr)^p &\text{if } -\frac12<r<\frac12. \end{cases} \eeq
In the last formula, $B_p$ is the complete Bell polynomial in $p$ variables (see \cite[Definition 2.4.1]{Peccati11}), $\zeta(z):= \sum_{k=1}^\infty k^{-z}$, for $z>1$, is the Riemann zeta function, and $\Ga(z):=\int_0^\infty y^{z-1}\ee^{-y}\,\dd y$ is Euler's gamma function.
\ethmA

The appearance of the Riemann zeta function in formula \eqref{const} is somewhat unexpected. It is related to the fact that in spatial dimension $d=1$, the $k$th eigenvalue of the Laplacian is proportional to $k^2$ and the norm on $H_r$ is defined using these eigenvalues (see \eqref{Hr}). Somewhat surprisingly, for $d\geq 1$ and a wide class of bounded open sets $D$,  conclusions similar to those of Theorem \ref{thmA} remain valid, after replacing the Riemann zeta function by the so-called spectral zeta function (see Corollaries \ref{cor1} and \ref{cor2} and Remarks \ref{rem2.7} and \ref{rem2.10}).

As we can see from both the formula for the normalizing sequence in \eqref{taunr0} and the formula for the limiting constant in \eqref{const}, the behavior of $V^{n,r}_p(u,t)$ changes at the critical value $r=-\frac12$. In fact, our analysis shows that for $r<-\frac12$, the terms in the sum \eqref{Vnp} become nearly independent and identically distributed as $n$ gets large, so that the convergence in \eqref{conv} can be interpreted as a \emph{(weak) law of large numbers}. By contrast, for $r\in[-\frac12, \frac12)$, they become \emph{deterministic} as $n\to\infty$, so the law of large numbers  is degenerate in this case.
The existence of a critical value of $r$ for the limiting behavior of $V^{n,r}_p(u,t)$ persists if we consider the more general equation \eqref{SHE}. In fact, Theorem~\ref{thmA} is a special case of our main Theorems~\ref{LLN} and \ref{LLN2} and their Corollaries~\ref{cor1} and \ref{cor2}.

Various authors have previously studied the  regularity of solutions to equations like \eqref{SHE}. For instance, the joint space-time Hölder regularity of the solution to the stochastic heat equation on $\R^d$ was established in \cite{SanzSole02}. Hölder regularity of the solution   as a process in $H_r$ (or related spaces) was investigated, for example,   in \cite[Theorem~11.8]{Peszat07} and \cite{vanNeerven08}. When $d=1$ and $D$ is an interval or $\R$, the power variations of $t\mapsto u(t,x)$, for fixed $x\in\R$, were analyzed  in \cite{Bibinger19b, Bibinger19, Cialenco20, Pospisil07,Swanson07} for $p\in\{2,4\}$ and in \cite{Chong20} for general powers. In the same setting, for fixed $t>0$, the power variations of $x\mapsto u(t,x)$ were discussed in \cite{Cialenco20,Hildebrandt19,Pospisil07} for $p=2$ and in \cite{Foondun15} for $\ga\in(\frac12,1]$ and $p=2/(2\ga-1)$. In the context of one-parameter stochastic processes, power variations have been investigated for semimartingales \cite{Jacod12}, fractional Brownian motion and related processes \cite{Corcuera06,Corcuera09, Nourdin08,Nourdin10, Nualart18}, and moving average processes \cite{BN11,Basse17, Corcuera13}, just to name a few.

To our best knowledge, power variations in Sobolev-type spaces have not been considered in the literature before. Even the existence of an $L^2$-continuous random field solution to \eqref{SHE} seems not yet to have been considered (our Proposition \ref{exist}), since it relies on fairly strong estimates concerning series of eigenfunctions of the Laplacian (see e.g.~Lemma \ref{A1}). Also, as Theorem~\ref{thmA} or the more general Theorems~\ref{LLN} and \ref{LLN2} below reveal, power variations in Sobolev-type spaces already show a rich behavior in first order. In particular, there is a phase transition in the law of large numbers that is not present in, for example, the power variations of fractional Brownian motion (cf.\ \cite{BN11,Corcuera06,Corcuera13}) or the power variations of the stochastic heat equation in time at a fixed spatial point (cf.\ \cite{Bibinger19,Chong20,Cialenco20,Swanson07}).

In the remainder of this paper, we write $A\lec B$ if there exists a finite constant $C>0$ (that is independent of all quantities of interest) such that $A\leq CB$. Moreover, we use the notation $\N=\{1,2,\dots\}$ and $\N_0=\{0,1,\dots\}$.

\section{Definitions and results}\label{sec2}

We first recall some basic definitions.
Let $0<\la_1\leq\la_2\leq \cdots$ and $(\phi_k)_{k\in\N}$ be the eigenvalues and corresponding normalized eigenfunctions of $-\Delta$, where $\Delta$ is the Dirichlet Laplacian on $D$ with zero boundary conditions, such that $(\phi_k)_{k\in\N}$ forms a complete orthonormal basis of $L^2(D)$, and each $\phi_k$ is smooth \cite[Corollary 8.11]{Gilbarg01} and bounded: see Lemma \ref{A1} (ii). As in 
\cite[Chapter~IV, Example~3]{Walsh86}), let $E_0$ be the set of $f$ of the form $f= \sum_{j=1}^N a_k \phi_k$ and define, for $r\in\R$,
$$
    \|f\|_{H_r} := \Biggl( \sum_{k=1}^{\infty} \la_k^r a_k^2 \Biggr)^{1/2}.
$$
Let $H_r:=H_r(D)$ be the completion of $E_0$ with respect to $\|\cdot \|_{H_r}$. Then $H_r$ is a Hilbert space, the \emph{Sobolev space of order $r$}, each element $\Phi$ of which can be identified with a series of the form
\beq\label{Hr}  \Phi=\sum_{k=1}^\infty a_k(\Phi)\phi_k, \qquad\text{where } a_k(\Phi) \in \R \text{ and } \|\Phi\|_{H_r} := \Biggl( \sum_{k=1}^{\infty} \la_k^r a_k(\Phi)^2 \Biggr)^{1/2}<\infty . \eeq
As noted in \cite[Remark~2.8]{Sheffield07}, the series defining $\Phi$ converges in the topology of the space of distributions on $D$ and in the $H_r$-norm. The inner product on $H_r$ is given by
\beq\label{Hr-ip} \langle \Phi_1,\Phi_2\rangle_{H_r}:=\sum_{k=1}^\infty \la_k^r a_k(\Phi)a_k(\Phi'),\qquad \Phi_1,\Phi_2\in H_r. \eeq
Moreover, $H_r\subseteq H_s$ for $s\leq r$, $H_0=L^2(D)$ with $\|\cdot\|_0=\|\cdot\|_{L^2(D)}$ by Parseval's identity, and the evaluation
$$
\langle \Phi_1, \Phi_2 \rangle := \sum_{k=1}^\infty a_k(\Phi_1) a_k(\Phi_2),\quad \Phi_1\in H_{-r},\quad \Phi_2\in H_r,
$$
puts $H_{-r}$ and $H_r$ in duality.

The \emph{spectral power of $-\Delta$ of order $\ga\in(0,\infty)$} is now defined   via
\beq\label{power} (-\Delta)^\ga\colon \bigcup_{r\in\R} H_r \to \bigcup_{r\in\R} H_r,\qquad (-\Delta)^\ga\Phi := \sum_{k=1}^\infty \la_k^\ga a_k(\Phi)\phi_k.  \eeq

It is standard to interpret the SPDE \eqref{SHE} using the notion of random field solution (see e.g.~\cite[Chapter 4]{Dalang20}), in which one gives a suitable meaning to the integral equation
\beq\label{eq:mild} u(t,x)  = \int_0^t \int_D g(t-s; x,y)\si(u(s,y))\, W( \dd s, \dd y)\qquad\text{a.s.,}  \eeq
which should be satisfied for all $(t,x)\in[0,\infty)\times D$, where  the Dirichlet Green's kernel $g$ takes the form
\begin{equation} \label{eq:def-heat}
g(t;x,y) :=   \sum_{k=1}^\infty \phi_k(x)\phi_k(y) \ee^{-\la^\ga_k t}\mathds 1_{t> 0},\qquad (t,x,y)\in [0,\infty)\times D^2.
\end{equation}

The existence and uniqueness of a mild solution to \eqref{SHE} can be established under appropriate conditions on $\ga$ and $D$. The proof of the next Proposition is given in the \hyperref[appn]{Appendix} after Lemma~\ref{A1}.

\bprop\label{exist} Let $D$ be a bounded open subset of $\R^d$ that satisfies the {\em cone property} (see \cite[Section~2]{Agmon65}). Let $\si$ be a globally Lipschitz function and let $\ga\in(\frac d2,\infty)$. Then there is a predictable random field $(t,x)\mapsto u(t,x)$, called the \emph{random field solution to \eqref{SHE}} that satisfies \eqref{eq:mild} for all $(t,x)\in(0,\infty)\times D$ and is such that
\beq\label{eq:Lp} \sup_{(t,x)\in[0,T]\times D} \E[|u(t,x)|^p] < \infty \eeq
for all $p,T>0$. In addition, $u$ is $L^p(\Om)$-continuous for all $p>0$.  Up to versions, $u$ is unique among all predictable random fields satisfying \eqref{eq:Lp} for $p=2$.
\eprop

\brem
A hyperbolic analogue of \eqref{SHE} (with $D=\R^d$ and a spatially homogeneous noise) was considered by the second author together with M.\ Sanz-Solé in \cite{Dalang05}.
\erem


In fact, our methods apply to a more general situation. Instead of studying \eqref{SHE} with a nonlinearity of a specific functional form, we shall consider
\beq\label{SHE-2}
\left\{\begin{array}{ll} \frac{\partial u}{\partial t}(t,x) = -(-\Delta)^\ga u(t,x)+\sigma(t,x) \dot W(t, x), & (t,x) \in [0,\infty)\times D,\\ u(t,x)=0, &  (t,x)\in[0,\infty)\times\partial D, \\
	u(0,x)=0, & x\in D, \end{array}\right.\eeq
with a general predictable random field $\si(t,x)$ satisfying mild regularity assumptions (so that \eqref{SHE} is included as a special case).
When $\ga >\frac d2$, the random field solution to \eqref{SHE-2} is given  by
\[
 u(t,x)= \int_0^t \int_D g(t-s;x,y) \si(s,y)\,W(\dd s,\dd y)
\]
for $t\in(0,\infty)$ and $x\in D$.
By the stochastic Fubini theorem \cite[Theorem 2.6]{Walsh86}, whose assumption is satisfied by \eqref{est1} because $\gamma > \tfrac{d}{2}$,
\beq
u(t,\cdot) =\sum_{k=1}^\infty  a_k(t)\phi_k \qquad\text{with}\qquad a_k(t) :=\int_0^t\int_D \ee^{-\la_k^\ga(t-s)}\phi_k(y)\si(s,y)\,W(\dd s,\dd y)
\label{akt}
\eeq
for $k\in\N$, where the series converges in $L^2(\Omega)$.

Another advantage of studying \eqref{SHE-2} is that we can drop the condition $\ga>\frac d2$ (which guarantees a random field solution and is needed for \eqref{SHE} in the case of multiplicative noise) if we allow $t\mapsto u(t,\cdot)$ to be a distribution-valued process. Indeed, by a direct calculation and Weyl's law (see \eqref{eq:Weyl}), we obtain for $t>0$ that
\beq\label{eq:L2finite}
\bbe\Biggl[\sum_{k=1}^\infty \la_k^r (a_k(t))^2 \Biggr] = \frac12\sum_{k=1}^\infty \la_k^{r-\ga} (1-\ee^{-2\la_k^\ga t})< \infty \iff r<\ga-\tfrac d2.
\eeq
Therefore, if $r<\ga-\tfrac d2$ and $\sigma$ satisfies \eqref{mom} below with $q=2$, then, for every $t\geq0$, the series in \eqref{akt} converges in $H_r$ almost surely and hence \emph{defines} a random variable $u(t,\cdot)$ in $H_r$ with $\E[\|u(t,\cdot)\|^2_{H_r}]< \infty$. As shown in \cite[Proposition~5.3]{Walsh86}, the $H_r$-valued process $t\mapsto u(t,\cdot)$ obtained in this way is  the \emph{weak solution} to \eqref{SHE-2} in the sense of \cite[Equation~(5.4)]{Walsh86} (see also \cite[Chapter~3]{Dalang20}). Let us remark that if $\ga=1$ and $\si\equiv 1$, then the weak solution to \eqref{SHE-2} is a dynamical analogue of the \emph{Gaussian free field} considered, for example, in \cite{Giacomin01, Sheffield07}.   The restriction  $r<\ga-\tfrac d2$ cannot be removed as the following result shows (see the \hyperref[appn]{Appendix} for a proof):
\bprop\label{uinHr} If $\si\equiv1$, then for all $t>0$, we have $u(t,\cdot)\in H_r$ a.s.~if and only if $r<\ga-\frac d2$.
\eprop

It is therefore natural to consider power variations in the spaces $H_r$ with $r< \ga-\tfrac d 2$. In fact,  more general functionals can be considered. Indeed, if $r<-\frac d2$, we shall investigate the \emph{(normalized) $F$-variations} of the solution $u$ to \eqref{SHE-2}, that is,
\beq\label{VnF} V^{n,r}_F(u,t):= \Delta_n\sum_{i=1}^{[t/\Delta_n]} F\biggl( \frac{u(i\Delta_n,\cdot)-u((i-1)\Delta_n,\cdot)}{\tau_n(r)} \biggr),\quad t\in[0,\infty),\quad n\geq1, \eeq
where $F\colon H_r\to\R$ is a continuous functional on $H_r$ satisfying certain regularity properties and the $\tau_n(r)$ are normalizing constants that we will define later in \eqref{taunr}.

If $-\frac d2\leq r<\ga-\frac d2$, we consider
\beq\label{Vnf} V^{n,r}_f(u,t):= \Delta_n\sum_{i=1}^{[t/\Delta_n]} f\biggl( \frac{\|u(i\Delta_n,\cdot)-u((i-1)\Delta_n,\cdot)\|_{H_r}}{\tau_n(r)} \biggr),\quad t\in[0,\infty),\quad n\geq1, \eeq
where $f\colon[0,\infty)\to\R$ is a continuous function satisfying certain regularity properties. For the reasons given in Remarks~\ref{rem-tight} and \ref{difference}, we do not consider functionals as general as \eqref{VnF} when $-\frac d2\leq r<\ga-\frac d2$.

For the case where $r<-\frac d2$, we   recall some facts about Gaussian measures on Hilbert spaces. Let $\call^{+,\mathrm{sym}}_1(H_r)$ be the space of all bounded linear operators $H_r\to H_r$ that are symmetric, nonnegative (meaning ``nonnegative definite'') and of finite trace (see \cite[Appendix~B]{Prevot07}). By \cite[Theorem~2.1.2]{Prevot07},  for every $Q\in \call^{+,\mathrm{sym}}_1(H_r)$, there is a Gaussian measure $N_r(0,Q)$ on $H_r$ with mean zero and covariance operator $Q$ (i.e., if $H\sim N_r(0,Q)$, then for every $h\in H_r$, the real-valued random variable $\langle H,h\rangle_{H_r}$ is normally distributed with mean $0$ and variance $\langle Qh,h\rangle_{H_r}$).

In the following, given a functional $F\colon H_r\to \R$ and some $Q\in\call^{+,\mathrm{sym}}_1(H_r)$, we write $\mu_F(Q) = \bbe[F(H)]$, where $H\sim N_r(0,Q)$, whenever this expectation exists. Still for $r<-\tfrac d2$, with a slight abuse of notation, given a nonnegative function $w\in L^\infty(D)$ (or a nonnegative random field $w$ with uniformly bounded expectation on $D$), we shall use the abbreviation $\mu_{r,F}(w)=\mu_F(Q_r(w))$ where the operator $Q_r(w)\colon H_r\to H_r$ is defined, for $h\in H_r$, by
\beq\label{Qphi1}
   Q_r(w) h := \sum_{k=1}^\infty \phi_k \sum_{\ell=1}^\infty \la_\ell^r a_\ell(h)\int_D \phi_k(y)\phi_\ell(y)w(y)\,\dd y,
\eeq
or, equivalently, for $h_1\in H_r$ and $h_2\in H_{-r}$, by
\beq\label{Qphi}
(Q_r(w) h_1)(h_2):= \sum_{k,\ell=1}^\infty  \la_\ell^r a_k(h_2)a_\ell(h_1)\int_D \phi_k(y)\phi_\ell(y)w(y)\,\dd y .\eeq
One easily checks that $Q_r(w) \in \call^{+,\mathrm{sym}}_1(H_r)$ since $r <-\tfrac d2$. A more tractable formula for the law $N_r(0,Q_r(w))$ is given in Lemma \ref{eNoQ}.

We first formulate our main result for $r < -\tfrac d2$. Given processes $(X^n(t))_{t\geq0}$ and $(X(t))_{t\geq0}$, we write $X^n\limq X$ (or $X^n(t)\limq X(t)$) if $\supt |X^n(t)-X(t)|\stackrel{L^q}{\longrightarrow}0$  for all $T>0$ as $n\to\infty$.

\bthm\label{LLN} Let $D$ be a bounded open subset of $\R^d$ with the cone property (see \cite[Section~2]{Agmon65}) and $\ga\in(0,\infty)$. Let $u$ be the solution to \eqref{SHE-2} given by \eqref{akt}. Suppose that $r<-\frac d2$ and let $\tau_n(r)=\sqrt{\Del}$.
Assume that $F\colon H_r \to\R$ is continuous with at most polynomial growth (i.e., there is $p>0$ such that $|F(h)|\lec 1+\|h\|_{H_r}^p$ for all $h\in H_r$). If $\si$ is a predictable random field that is continuous in probability and satisfies
\beq\label{mom} \sup_{(t,x)\in[0,T]\times D} \bbe[|\si(t,x)|^q] < \infty
 \qquad\text{for all } q>0\text{ and } T>0,
\eeq
then for every $q>0$, as $n\to\infty$, $V^r_F(u,t) := L^q\text{-}\lim_{n\to\infty} V^{n,r}_F(u,t)$ exists,
\beq\label{eq:LLN} V^{n,r}_F(u,t) \limq V^r_F(u,t)\qquad\text{and}\qquad V^r_F(u,t)= \int_0^t \mu_{r,F} \Bigl( \si^2(s,\cdot)\Bigr)\,\dd s. \eeq
\ethm
Theorem \ref{LLN} is proved in Section \ref{sec3.1}.

\brem\label{L2-rem} By well-known results on uniform integrability (combine Proposition~4.12 with the discussion before Lemma~4.10 in \cite{Kallenberg02}), if $\si$ is continuous in probability and satisfies \eqref{mom}, then it is automatically $L^p(\Om)$-continuous for every $p\geq1$. In particular, for every $p\geq1$ and $T>0$,
\beq\label{L2-cont} w_p(\eps;T):= \sup\Bigl\{(\bbe[|\si(t,x)-\si(s,y)|^p])^{\frac1p}: |t-s|+|x-y|<\eps,~s,t\in[0,T]\Bigr\} \to 0
\eeq
as $\eps\to0$.
\erem

Let us apply Theorem~\ref{LLN} to the power variations in \eqref{Vnp} of \emph{even} order, which correspond to the functionals $F_{2p}(h) = \|h\|_{H_r}^{2p}$ for some $p\in\N$. In this case, $\mu_{F_{2p}}$ can be determined explicitly. To this end, we introduce
for a symmetric matrix $A\in\R^{p\times p}$, its \emph{$\al$-permanent} (see \cite{VereJones97}),  defined as
\beq\label{perm} \mathrm{per}_\al(A):=\sum_{\Sigma\in\mathrm{Sym}_p} \al^{\#\Sigma} \prod_{i=1}^p A_{i\Sigma(i)},\qquad \al\in\R, \eeq
where $\mathrm{Sym}_p$ is the symmetric group of all permutations $\Sigma$ on $\{1,\ldots,p\}$ and $\#\Sigma$ is the number of cycles in $\Sigma$. If $\al=1$, this is the ordinary permanent of $A$, while for $\al=-1$, this is $(-1)^n$ times the determinant of $A$. For $\al=\frac12$, this matrix function has an important connection to moments of Gaussian random variables, see  \cite[Equation~(3.23)]{McCullagh18}: If $X_1,\ldots,X_p$ are jointly Gaussian with mean $0$ and covariance matrix $C$, then
\beq\label{prodmom} \bbe[X_1^2\cdots X_p^2] = 2^p\mathrm{per}_{\frac12}(C). \eeq
In the following, we will apply the permanent function to matrices of the form
\[ D_r(y_1,\ldots,y_p) := \begin{pmatrix} \langle \delta_{y_1},\delta_{y_1}\rangle_{H_r} &\dots& \langle \delta_{y_1},\delta_{y_p}\rangle_{H_r} \\ \vdots & \ddots & \vdots \\ \langle \delta_{y_p},\delta_{y_1}\rangle_{H_r} &\dots& \langle \delta_{y_p},\delta_{y_p}\rangle_{H_r}  \end{pmatrix}\]
and define
\beq\label{perm2} \mathrm{per}^{(r)}_\al(y_1,\ldots,y_p) := \mathrm{per}_{\al}(D_r(y_1,\ldots,y_p)), \eeq
where $r<-\frac d2$, $y_1,\dots,y_p\in D$ are given spatial points and $\delta_y$ is the Dirac delta function at $y$. Notice that by \eqref{est1}, $\delta_{y} \in H_r$ for $r < - \frac d2$, and
\[
   \langle \delta_{y_i},\delta_{y_j}\rangle_{H_r} =  \sum_{k=1}^\infty \la_k^r \phi_k(y_i)\phi_k(y_j).
\]

\bcor\label{cor1} Assume the hypotheses of Theorem~\ref{LLN} (in particular, $r<-\frac d2$).
\benu
\item The process $t\mapsto u(t,\cdot)$ has a locally finite quadratic variation in $H_r$: for any $q>0$,
\beq\label{2var}\begin{aligned} V^{n,r}_2(u,t)&=\sum_{i=1}^{[t/\Del]} \|u(i\Delta_n,\cdot)-u((i-1)\Delta_n,\cdot)\|^2_{H_r}\\ & \limq V^r_2(u,t)=\int_0^t \int_D \si^2(s,y)\|\delta_y\|_{H_r}^2\,\dd y\,\dd s. \end{aligned}  \eeq
\item More generally, if $p\in\N$, then for any $q>0$, $ V^{n,r}_{2p}(u,t) \limq  V^r_{2p}(u,t)$, where
$$
V^{n,r}_{2p}(u,t)=\Del^{1-p} \sum_{i=1}^{[t/\Del]} \|u(i\Delta_n,\cdot)-u((i-1)\Delta_n,\cdot)\|^{2p}_{H_r}
$$
and
\beq\label{2pvar}
\begin{split}
 V^r_{2p}(u,t)&=\int_0^t \sum_{\Sigma\in\mathrm{Sym}_p} 2^{p-\#\Sigma} \int_D \cdots\int_D \langle\delta_{y_1},\delta_{y_{\Sigma^{-1}(1)}}\rangle_{H_r}\cdots\langle\delta_{y_p},\delta_{y_{\Sigma^{-1}(p)}}\rangle_{H_r}\\
	&\qquad\times \si^2(s,y_1)\cdots\si^2(s,y_p)\,\dd y_1\cdots\,\dd y_p\,\dd s\\
	&=2^p\int_0^t \int_D \cdots\int_D \mathrm{per}^{(r)}_{\frac12}(y_1,\ldots,y_p)\,  \si^2(s,y_1)\cdots\si^2(s,y_p)\,\dd y_1\cdots\,\dd y_p \,\dd s.
	\end{split}
\eeq
\item In particular, in the case of additive noise where $\si(s,y)\equiv \si$ for some $\si\in\R$, we have
\beq\label{2pvar-add}
 V^r_{2p}(u,t)=\si^{2p} 2^p B_p(x_1,\dots,x_p)t,
\eeq
where $B_p$ is the complete Bell polynomial in $p$ variables (see \cite[Definition 2.4.1]{Peccati11}),
\beq\label{xl}
   x_\ell := \frac{(\ell-1)!}{2} \zeta_D(-\ell r),\qquad \ell=1,\dots,p,
\eeq
and $\zeta_D(z):= \sum_{k=1}^\infty \la_k^{-z}$. In particular, for $p=1$, the quadratic variation in $H_r$ of $t\mapsto u(t,\cdot)$ is $V^r_{2}(u,t) = \si^2 \zeta_D(-r) t$.
\eenu
\ecor

Corollary \ref{cor1} is proved at the end of Section \ref{sec3.1}.

\brem\label{rem2.7} The function $\zeta_D(z):= \sum_{k=1}^\infty \la_k^{-z}$ is known as the  {\em spectral zeta function} or as the \emph{Minakshisundaram--Pleijel zeta function} in the mathematical physics literature (see \cite{Elizalde12,Hawking77,Voros87,Voros92}). By Weyl's law (see \eqref{eq:Weyl}), $\zeta_D(z)$ well defined when $\mathrm{Re}(z)>\frac d2$. In the special case where $d=1$ and $D=[0,\pi]$, then $\la_k = k^2$ and $\zeta_D(z) = \zeta(2z)$, where $\zeta(\cdot)$ is the classical Riemann zeta function. When $\si(s,y) \equiv \sigma$, the formula $V^r_{2}(u,t) = \si^2 \zeta_D(-r) t$ shows that knowing the  quadratic variation of $u$ in $H_r$ for $r<-\tfrac{d}{2}$ yields, via the function $\zeta_D$, information about the eigenvalues $(\la_k)_{k\in\N}$ and hence about the domain $D$.
\erem

Next, in the case where $-\frac d2\leq r<\ga-\frac d2$, we shall only consider functionals of the form $F(h)=f(\|h\|_{H_r})$ with a real-valued function $f\colon [0,\infty)\to\R$; see \eqref{Vnf}.

\bthm\label{LLN2} Let $\ga\in(0,\infty)$ and $D$ be a bounded connected open subset of $\R^d$ with a piecewise smooth boundary in the sense of \cite[Definition~1.17]{Hezari18} and the cone property (see \cite[Section~2]{Agmon65}). Suppose that $-\frac d2\leq r<\ga-\frac d2$ and let
\beq\label{taunr} \tau_n(r):= \begin{cases} \Delta_n^{\frac1{2\ga}(\ga-\frac{d}{2}-r)}&\text{if } -\frac d2<r<\ga-\frac d2,\\ (\Del\lvert\log\Del\rvert)^{\frac12}&\text{if } r=-\frac d2.\end{cases}\eeq
Assume that $f\colon[0,\infty)\to\R$ is continuous with at most polynomial growth (i.e., there is $p>0$ such that $|f(x)|\lec 1+x^p$ for all $x\geq0$).

Let $u$ be the solution to \eqref{SHE-2} given by \eqref{akt}, and assume that $\si$ is a predictable random field that is continuous in probability and satisfies \eqref{mom}.
Then, as $n\to\infty$, $V^r_f(u,t) := L^q\text{-}\lim_{n\to\infty} V^{n,r}_f(u,t)$ exists, 
\beq\label{eq:LLN-2} V^{n,r}_f(u,t) \limq V^r_f(u,t)\qquad\text{and}\qquad V^r_f(u,t)= \int_0^t f\Biggl( \sqrt{\frac{K_r}{|D|}\int_D \si^2(s,y)\,\dd y}\Biggr)\,\dd s \eeq
 for all $q>0$,
where
\beq\label{Kr} K_r := 
\begin{cases} \frac{|D|}{(4\pi)^{d/2} \Ga(\frac d2)(\ga-\frac{d}{2}-r)} \Ga\Bigl(\frac1\ga( r  +\frac{d}{2})\Bigr)
	&\text{if } -\frac d2<r<\ga-\frac d2,\\ \frac{|D|}{\ga(4\pi)^{d/2} \Ga(\frac d2)} &\text{if } r=-\frac d2,\end{cases} 
\eeq
and where $|D|$ denotes the Lebesgue measure of the set $D$. 
\ethm

Theorem \ref{LLN2} is proved in Section \ref{nsec3.3}. Using the functions $f(x) = x^p$, $p>0$, we immediately get the following corollary.

\bcor\label{cor2} Assume the hypotheses of Theorem~\ref{LLN2} and let $-\frac d2< r<\ga-\frac d2$.
\benu
	\item The process $t\mapsto u(t,\cdot)$ has a locally finite $\frac{2\ga}{\ga- d/2 -r}$-variation in $H_r$. More precisely, for any $q>0$,
	\beq\label{exact-var}\begin{aligned} V^{n,r}_{\frac{2\ga}{\ga- d/2 -r}}(u,t)&=\sum_{i=1}^{[t/\Del]} \|u(i\Delta_n,\cdot)-u((i-1)\Delta_n,\cdot)\|^{\frac{2\ga}{\ga- d/2 -r}}_{H_r}\\ & \limq V^r_{\frac{2\ga}{\ga- d/2 -r}}(u,t)=\biggl( \frac{K_r}{|D|}\biggr)^{\frac{\ga}{\ga- d/2 -r}}\int_0^t \biggl( \int_D \si^2(s,y)\,\dd y\biggr)^{\frac{\ga}{\ga- d/2 -r}}\,\dd s. \end{aligned}  \eeq
	\item More generally, for any $p>0$,
	\beq\label{p-var}\begin{aligned} V^{n,r}_p(u,t)&=\Del^{1-\frac p{2\ga}(\ga-\frac{d}{2}-r)}\sum_{i=1}^{[t/\Del]} \|u(i\Delta_n,\cdot)-u((i-1)\Delta_n,\cdot)\|^p_{H_r}\\ & \limq V^r_p(u,t)=\biggl( \frac{K_r}{|D|}\biggr)^{\frac p2}\int_0^t \biggl( \int_D\si^2(s,y)\,\dd y\biggr)^{\frac p2}\,\dd s. \end{aligned}  \eeq
	\item In particular, in the case of additive noise where $\si(s,y)\equiv \si$ for some $\si\in\R$, we have
		\beq\label{p-var-add}V^r_p(u,t)=|\si|^p K_r^{\frac p2}t. \eeq
	\item Part (ii) and part (iii) remain valid for $r=-\frac d2$ when we replace the first line in \eqref{p-var} by
	\beq\label{p-var-12} V^{n,r}_p(u,t)=\frac{\Del^{1-\frac{p}{2}}}{\lvert\log \Del\rvert^{\frac p2}}\sum_{i=1}^{[t/\Del]} \|u(i\Delta_n,\cdot)-u((i-1)\Delta_n,\cdot)\|^p_{H_r}.  \eeq
\eenu
\ecor

\brem\label{rem2.10} Contrary to the case $r < -\tfrac{d}{2}$, in the case where $-\tfrac{d}{2} \leq r < \gamma - \tfrac{d}{2}$, the quadratic variation of $u$ is infinite and the $\frac{2\ga}{\ga- d/2 -r}$-variation, which is finite, contains very little information about the domain $D$, other than its volume $|D|$, via the constant $K_r$.
\erem
Our proofs also yield the optimal Hölder regularity of the sample paths of $u$ (as a process with values in $H_r$); cf.\ also \cite[Theorem~11.8]{Peszat07} and \cite{vanNeerven08}. Similar results in the hyperbolic case (with $D=\R^d$) were obtained in \cite{Dalang05}. The proof of the following Corollary \ref{Holder} is given at the end of Section~\ref{nsec3.4}.
\bcor\label{Holder} Let $r<\ga-\frac d2$ and
assume the hypotheses of Theorem~\ref{LLN} (if $r<-\frac d2$) or Theorem~\ref{LLN2} (if $r\geq -\frac d2$). Then $t\mapsto u(t,\cdot)$, viewed as a stochastic process with values in $H_r$ has a version that is almost surely locally $\al$-Hölder continuous for all $\al<\al(r)$,
where $\al(r):=\frac12$ if $ r\leq -\frac d2$ and $\al(r):=\frac{1}{2\ga}(\ga-\frac d2 -r)$ if $-\frac d2<r<\ga-\frac d2$.
Moreover, if $\si\equiv1$ and $\al>\al(r)$, then $t\mapsto u(t,\cdot)$ does not have a version that is almost surely locally $\al$-Hölder continuous.
\ecor


The next section contains the proofs of all the results that we have stated above. In Section~\ref{compare}, we also explain in Remarks \ref{rem-tight} and \ref{difference} why in Theorem~\ref{LLN2} (where $-\frac d2\leq r<\ga-\frac d2$), we did not study functionals as general as in Theorem~\ref{LLN} (where $r<-\frac d2$). In the \hyperref[appn]{Appendix}, we gather several important properties and estimates concerning eigenvalues and eigenfunctions of the Laplacian on $D$, which are used in these proofs.

\section{Proofs}

We first show that if $\si\equiv 1$, then the $H_r$-norm of an increment $u(i\Delta_n,\cdot)-u((i-1)\Delta_n,\cdot)$ is typically of order $\tau_n(r)$ as defined in Theorem \ref{LLN} or \eqref{taunr} (depending on the value of $r$).
\blem\label{size-incr} Let
\beq\label{Kr-2} K_r := \sum_{k=1}^\infty \la_k^r = \zeta_D(-r), \qquad\text{if } r<-\frac d2, \eeq
and let $K_r$ be defined as in \eqref{Kr} if $-\frac d2\leq r<\ga-\frac d2$. If $u$ is the solution to \eqref{SHE-2} with $\si\equiv1$, then for every $\eps>0$,
\beq\label{limit} \lim_{n\to\infty} \sup_{i\colon i\Del\geq\eps} \Biggl\lvert\frac{\bbe[\|u(i\Delta_n,\cdot)-u((i-1)\Delta_n,\cdot)\|_{H_r}^2]}{\tau_n(r)^2} -K_r\Biggr\rvert = 0.\eeq
\elem
\bpr  By standard calculations, and using the notation $\iint_a^b = \int_a^b \int_D$,
\beq\label{calc}
\begin{aligned}
&\bbe[\|u(i\Delta_n,\cdot)-u((i-1)\Delta_n,\cdot)\|_{H_r}^2]\\ 
&\quad= \sum_{k=1}^\infty \la_k^{r}\bbe[(a_k(i\Del)-a_k((i-1)\Del))^2] \\
&\quad=\sum_{k=1}^\infty \la_k^{r} \Biggl( \iint_0^{i\Del} \ee^{-2\la_k^\ga (i\Del-s)}\phi_k(y)^2\,\dd s\,\dd y+\iint_0^{(i-1)\Del} \ee^{-2\la_k^\ga((i-1)\Del-s)}\phi_k(y)^2\,\dd s\,\dd y\\
&\quad\qquad-2\iint_0^{(i-1)\Del} \ee^{-\la_k^\ga(i\Del-s)-\la_k^\ga((i-1)\Del-s)}\phi_k(y)^2\,\dd s\,\dd y \Biggr)\\
&\quad=\sum_{k=1}^\infty \la_k^{r}\Biggl( \int_0^{i\Del} \ee^{-2\la_k^\ga s}\,\dd s -2\int_0^{(i-1)\Del} \ee^{-\la_k^\ga (s+\Del)-\la_k^\ga s}\,\dd s+\int_0^{(i-1)\Del} \ee^{-2\la_k^\ga s}\,\dd s   \Biggr)\\
&
\quad= \displaystyle \sum_{k=1}^\infty \la_k^{r} \frac{1-\ee^{-2\la_k^\ga i\Del}-2\ee^{-\la_k^\ga \Del}+2\ee^{-\la_k^\ga \Del}\ee^{-2\la_k^\ga (i-1)\Del}+1-\ee^{-2\la_k^\ga (i-1)\Del}}{2\la_k^\ga}\\
&\quad=\sum_{k=1}^\infty \la_k^{r-\ga}(1-\ee^{-\la_k^\ga\Del}) -\frac12 \sum_{k=1}^\infty \la_k^{r-\ga} \ee^{-2\la_k^\ga i\Del}  (\ee^{\la_k^\ga\Del}-1)^2.
\end{aligned}\eeq
Let us denote the first and the second series by $A^{n,1}_r$ and $A^{n,i,2}_r$, respectively. By dominated convergence, we have
\beq\label{second-term} \limsup_{n\to\infty} \sup_{i\colon i\Del\geq\eps} \Delta_n^{-2}A^{n,i,2}_r \leq \lim_{n\to\infty}  \sum_{k=1}^\infty \frac{\la_k^{r-\ga}}{2} \ee^{-2\la_k^\ga \eps}  \Biggl(\frac{\ee^{\la_k^\ga\Del}-1}{\Del} \Biggr)^2 = \sum_{k=1}^\infty \frac{\la_k^{r+\ga}}{2}\ee^{-2\la_k^\ga \eps}<\infty, \eeq
which shows that the contribution of $A^{n,i,2}_r$ is negligible for the limit taken in \eqref{limit}.

Concerning $A^{n,1}_r$, if $r<-\frac d2$, then we immediately deduce from dominated convergence that
\beq\label{convKr} \lim_{n\to\infty} \Del^{-1}A^{n,1}_r = \sum_{k=1}^\infty \la_k^r = K_r,  \eeq
which implies \eqref{limit}.
 If $-\frac d2\leq r<\ga-\frac d2$, then Lemma~\ref{convdens} (i) shows $\tau(n)^{-2}A^{n,1}_r=H_n\to K_r$.
\epr

\subsection{Proof of Theorem~\ref{LLN} and Corollary~\ref{cor1} ($r < -\tfrac d2$)}\label{sec3.1}

\bpr[Proof of Theorem~\ref{LLN}]
We only show \eqref{eq:LLN} for $q=1$. By \cite[Proposition~4.12]{Kallenberg02} and the discussion before Lemma~4.10 in \cite{Kallenberg02}, the statement for all larger values of $q$ follows easily from the hypothesis that $F$ has polynomial growth and the fact that $u$ has locally uniformly bounded moments of all orders, by using the fact that convergence in probability together with bounded $L^{q+\eps}$-moments implies convergence in $L^q$.

Since $V^{n,r}_F$ and $V^r_F$ in \eqref{eq:LLN} are linear in $F$, by decomposing $F$ into its positive and negative parts, we may further assume without loss of generality that $F$ is nonnegative. In this case, both $V^{n,r}_F(u,t)$ and $V^r_F(u,t)$ are nonnegative and increasing in $t$, so uniform convergence on compacts is equivalent to pointwise convergence in $t$ (see \cite[Chapter VI, Theorem 2.15 c)]{Jacod03}). Hence, in the following, we consider a fixed time point $t>0$. Furthermore, for brevity, we use the notation
\beq\label{abbr}\begin{aligned}
	\Delta^n_i u &:= u(i\Delta_n,\cdot)-u((i-1)\Delta_n,\cdot),\quad \Delta^n_i a_k := a_k(i\Del)-a_k((i-1)\Del),\\
	\Delta^n_i g(s;x,y)&:=g(i\Del-s;x,y)-g((i-1)\Del-s;x,y)\bone_{[0,(i-1)\Del]}(s),\\
	\Delta^n_i \ee_k(s)&:=\ee^{-\la_k^\ga(i\Del-s)}\bone_{[0,i\Del]}(s)-\ee^{-\la_k^\ga((i-1)\Del-s)}\bone_{[0,(i-1)\Del]}(s).\end{aligned}
\eeq

Recall from \eqref{akt} that the coefficients $a_k(t)$ of $u(t,\cdot)$ are in fact Ornstein--Uhlenbeck-type processes,
that is, for every $k\in\N$, $a_k$ satisfies the stochastic differential equation
\beq\label{SDE} \dd a_k(t)= -\la_k^\ga a_k(t)\,\dd t + \dd X_k(t),\qquad a_k(0)=0,\eeq
where
\beq\label{Xkt} X_k(t):=\int_0^t\int_D \phi_k(y)\si(s,y)\,W(\dd s,\dd y).\eeq
%
As a first step, we will show that the drift part is asymptotically negligible, leading us to approximate $u(t,\cdot)$ by $u'(t,\cdot)$, where
\beq\label{uprime}
   u'(t,\cdot):= \sum_{k=1}^\infty \phi_k X_k(t),\qquad t\in[0,\infty),
   \eeq
and the limit in \eqref{uprime} is taken in $H_r$.

As a second step, we locally freeze the coefficient $\sigma$, leading us to approximate $V^{n,r}_F(u',t)$ by
$$
\Del\sumt F\biggl( \frac{h^n_i}{\sqrt{\Del}}\biggr),
$$
where, for $n\in\N$ and $i=1,\ldots,[t/\Del]$, $h^n_i$ is the $H_r$-valued random variable
\beq\label{hni} h^n_i := \sum_{k=1}^\infty \phi_k\iint_{(i-1)\Del}^{i\Del}\phi_k(y)\si((i-1)\Del,y)\,W(\dd s,\dd y).  \eeq
With this notation, we will prove that \eqref{eq:LLN} holds by writing
\beq\label{eq3.33}
   V^{n,r}_F(u,t) - \int_0^t \mu_F( \si^2(s,\cdot))\,\dd s = I_{n,1} + \cdots + I_{n,4},
\eeq
where
\begin{align*}
   I_{n,1} &:= V^{n,r}_F(u,t)-V^{n,r}_F(u',t),
  \qquad I_{n,2} := V^{n,r}_F(u',t) - \Del\sumt F\biggl( \frac{h^n_i}{\sqrt{\Del}}\biggr),\\
   I_{n,3} &:= \Del\sumt F\biggl( \frac{h^n_i}{\sqrt{\Del}}\biggr) - \Del\sumt \mu_F(Q^{n,i}_r),\\
  I_{n,4} &:= \Del\sumt \mu_F(Q^{n,i}_r) - \int_0^t \mu_F( \si^2(s,\cdot))\,\dd s,
\end{align*}
and the operator $Q^{n,i}_r$ is defined in \eqref{Qni} below. The terms $I_{n,1}$ and $I_{n,2}$ are exactly the errors incurred by the two approximations described above. Their convergence to $0$ in $L^1(\Om)$ will be proved in Lemmas \ref{remove-drift} and \ref{discr-sigma}, respectively, based on a tightness argument (see Lemma~\ref{tight}). Next, we establish the law of large numbers for $\Del\sumt F( {h^n_i}/{\sqrt{\Del}})$ in Lemma~\ref{LLN-core-1}, which, together with a calculation of the conditional expectation that appears there (see Lemma~\ref{comp-condexp}), shows that $\E[|I_{n,3}|]\to0$ as $n\to\infty$. Finally, Lemma~\ref{cex-conv} shows that the discrete sum that we obtain from the previous step approximates $V^r_F(u,t)=\int_0^t\mu_F(\si^2(s,\cdot))\,\dd s$ in $L^1(\Omega)$, that is, $\E[|I_{n,4}|]\to 0$. To summarize, once Lemmas~\ref{remove-drift}--\ref{cex-conv} are proved, the proof of Theorem \ref{LLN} will be complete.
\epr

\blem\label{eNoQ}
Let $r < -\tfrac d2$. If $H$ is a random vector with values in $H_r$ and law $N_r(0,Q_r(w))$, where $Q_r(w)$ is defined in \eqref{Qphi1} and \eqref{Qphi}, then $H$ has the same law as $\sum_{k=1}^\infty X_k\, \la_k^{- r/2} \phi_k$, where the $X_k$ are jointly Gaussian centered random variables with covariances
\beq\label{covkl}
\cov(X_k,X_\ell)  =  \la_k^{\frac r2}\la_\ell^{\frac r2} \int_D \phi_k(y)\phi_\ell(y)w(y)\,\dd y.
\eeq
\elem
\bpr
Consider the orthonormal basis of $H_r$ given by
\beq\label{basis_b}
\Bigl(b_k:=\la_k^{-\frac r2}\phi_k: k\in\N\Bigr)
\eeq
 and note that $H =\sum_{k=1}^\infty a_k(H) \la_k^{r/2}b_k$. Define $X_k = \la_k^{r/2}a_k(H)$. Since $a_k(b_\ell) = \la_k^{-r/2}\delta_{k,\ell}$, we have $X_k = \langle H,b_k\rangle_{H_r}$, and by \cite[Theorem 2.1.2]{Prevot07},
\[
\cov(X_k,X_\ell) = \E[\langle H,b_k\rangle_{H_r}\, \langle H,b_\ell\rangle_{H_r}]= \langle Q_r(w)b_k, b_\ell \rangle_{H_r}.
\]
By \eqref{Qphi1},
\begin{align*}
Q_r(w)b_k &=  \sum_{k_1=1}^\infty \phi_{k_1}\, \la_k^r a_k(b_k)\int_D \phi_{k_1}(y)\phi_k(y)w(y)\,\dd y\\
&=  \sum_{k_1=1}^\infty b_{k_1} \, \la_{k_1}^{\frac r2} \, \la_k^{\frac r2} \int_D \phi_{k_1}(y)\phi_k(y)w(y)\,\dd y.
\end{align*}
Therefore,
\begin{align*}
\langle Q_r(w)b_k, b_\ell \rangle_{H_r} 
& =   \la_\ell^{\frac r2}  \la_k^{\frac r2}   \int_D \phi_{\ell}(y)\phi_k(y)w(y)\,\dd y.
\end{align*}
This completes the proof.
\epr

\bpr[Proof of Corollary~\ref{cor1}] We first prove part (ii). The $2p$-variation obviously corresponds to the functional $F\colon H_r\to \R$ given by $F(h)=\|h\|_{H_r}^{2p}$, which satisfies the hypotheses of Theorem~\ref{LLN}. For a deterministic function $w \in L^2(D)$, let $H$ be a random vector in $H_r$ with law $N_r(0,Q_r(w))$. By definition and Lemma \ref{eNoQ},
\beq\label{mu2p}
\begin{aligned}\mu_{r,F}(w) &= \bbe[F(H)] = \bbe[\|H\|_{H_r}^{2p}] = \bbe\Biggl[ \Biggl(\sum_{k=1}^\infty X_k^2 \Biggl)^p \Biggr] 
	=\sum_{k_1,\ldots,k_p=1}^\infty \bbe\Big[ X_{k_1}^2\cdots X_{k_p}^2\Big],
\end{aligned}
\eeq
where 
$(X_k)_{k\in\N}$ is a sequence of jointly Gaussian centered random variables with covariances  given by \eqref{covkl}.
Using \eqref{prodmom} to evaluate the expectations in the last line of \eqref{mu2p}, we obtain
\begin{align*}
\mu_{r,F}(w) &= 2^p\sum_{k_1,\ldots,k_p=1}^\infty \mathrm{per}_{\frac12}\Bigl(\cov(X_{k_i},X_{k_j})_{i,j=1}^p\Bigr)\\
&=\sum_{\Sigma\in\mathrm{Sym}_p} 2^{p-\#\Sigma} \int_D \cdots\int_D \sum_{k_1,\ldots,k_p=1}^\infty \la_{k_1}^{\frac r2}\, \la^{\frac r 2}_{k_{\Sigma(1)}}\cdots \la_{k_p}^{\frac r2} \, \la_{k_{\Sigma(p)}}^{\frac r2}\\
&\qquad\times \phi_{k_1}(y_1)\phi_{k_{\Sigma(1)}}(y_1)\cdots \phi_{k_p}(y_p) \phi_{k_{\Sigma(p)}}(y_p)  w(y_1)\cdots w(y_p)\,\dd y_1\cdots\,\dd y_p.
\end{align*}
Since $\Sigma$ is one-to-one, each factor $\la_{k_i}^{r/2}$, for $i=1,\ldots,p$, appears exactly twice, and the same holds for the function $\phi_{k_i}$, but once with argument $y_i$ and once with argument $y_{\Sigma^{-1}(i)}$. Thus,
\beq\label{muFsigma}\begin{aligned}
	\mu_{r,F}(\si^2(s,\cdot)) &= \sum_{\Sigma\in\mathrm{Sym}_p} 2^{p-\#\Sigma} \int_D \cdots\int_D \sum_{k_1=1}^\infty \la_{k_1}^{r} \phi_{k_1}(y_1)\phi_{k_1}(y_{\Sigma^{-1}(1)}) \\
	&\qquad\times\cdots\times \sum_{k_p=1}^\infty \la_{k_p}^{r} \phi_{k_p}(y_p)\phi_{k_p}(y_{\Sigma^{-1}(p)})  \si^2(s,y_1)\cdots\si^2(s,y_p)\,\dd y_1\cdots\,\dd y_p\\
	&=\sum_{\Sigma\in\mathrm{Sym}_p} 2^{p-\#\Sigma} \int_D \cdots\int_D \langle\delta_{y_1},\delta_{y_{\Sigma^{-1}(1)}}\rangle_{H_r}\cdots\langle\delta_{y_p},\delta_{y_{\Sigma^{-1}(p)}}\rangle_{H_r}\\
	&\qquad\times \si^2(s,y_1)\cdots\si^2(s,y_p)\,\dd y_1\cdots\,\dd y_p,
\end{aligned}\eeq
and part (ii) of Corollary~\ref{cor1} follows from Theorem~\ref{LLN}.

The other two parts of the corollary are special cases of the result that we have just proved. Indeed, for $p=1$, the permanent of a number $A\in\R$ is simply $\mathrm{per}_\al(A)=\al A$. And because $D_r(y)=\|\delta_y\|_{H_r}^2$, part (i) follows.

Regarding part (iii), if $\si(s,y)\equiv \si$, then \eqref{muFsigma} becomes
\beq\label{3.37}\begin{aligned}
	\mu_{r,F}(\si^2) &= \si^{2p}\sum_{\Sigma\in\mathrm{Sym}_p} 2^{p-\#\Sigma} \int_D \dd y_1 \cdots\int_D \dd y_p \Biggl[\sum_{k_1=1}^\infty \la_{k_1}^{r} \phi_{k_1}(y_1)\phi_{k_1}(y_{\Sigma^{-1}(1)}) \Biggr] \\
	&\qquad \times\cdots\times \Biggl[\sum_{k_p=1}^\infty \la_{k_p}^{r} \phi_{k_p}(y_p)\phi_{k_p}(y_{\Sigma^{-1}(p)}) \Biggr] \\
	&= \si^{2p}\sum_{\Sigma\in\mathrm{Sym}_p} 2^{p-\#\Sigma} \sum_{k_1=1}^\infty \cdots \sum_{k_p=1}^\infty  \la_{k_1}^{r} \cdots \la_{k_p}^{r}  \int_D \dd y_1 \cdots\int_D \dd y_p \, \phi_{k_1}(y_1)\phi_{k_1}(y_{\Sigma^{-1}(1)}) \\
	&\qquad  \times \cdots \times \phi_{k_p}(y_p)\phi_{k_p}(y_{\Sigma^{-1}(p)}).
\end{aligned}\eeq

Suppose that $\#\Sigma = \ell$ and these $\ell$ cycles have respective lengths $q_1,\dots, q_\ell$, with $q_1 + \cdots + q_\ell = p$. From the orthogonality property of the functions $\phi_{k_i}$, we see that the integral in \eqref{3.37} will vanish unless
\[
k_1 = k_{\Sigma(1)} = k_{\Sigma(\Sigma(1))} = \cdots,
\]
that is, the indices $k_i$ in each cycle of $\Sigma$ coincide.

The term in \eqref{3.37} corresponding to such a $\Sigma$ is
\beq\label{3.38}
2^{p-\ell} \Biggl[\sum_{i_1=1}^\infty \la_{i_1}^{r q_1} \Biggr] \cdots \Biggl[\sum_{i_\ell=1}^\infty \la_{i_\ell}^{r q_\ell} \Biggr].
\eeq
We would like to determine the number of permutations $\Sigma$ in $\mathrm{Sym}_p$ which consist of $\ell$ cycles  with respective lengths $q_1,\dots, q_\ell$. In order to obtain such a permutation of $\{1,\dots,p\}$, we first build, using the notation of \cite[Section~2.1]{Peccati11}, a partition of $\{1,\dots,p\}$ with $r_i$ blocks of size $i$ for every $i=1,\dots,p$, where $r_1,\dots,r_p \in \N_0$ are such that $1 r_1 + 2 r_2 + \cdots + p r_p = p$. According to  \cite[(2.3.8)]{Peccati11}, the number of such partitions is
\[
\Biggl[\begin{array}{c} p \\ r_1,\dots,r_p \end{array} \Biggr] := \frac{p!}{(1!)^{r_1}\, r_1! \, (2!)^{r_2}\, r_2!\,  \cdots\, (p!)^{r_p}\, r_p!}.
\]
For each of the $r_i$ blocks of size $i$, there are $(i-1)!$ possible cycles formed with the elements of this block. For each choice of these cycles, we obtain a permutation $\Sigma \in \mathrm{Sym}_p$ consisting of $r_1$ cycles of length $1,\dots, r_p$ cycles of length $p$. For this permutation, the term \eqref{3.38} is equal to
\begin{align*}
& 2^{p-(r_1+ \cdots + r_p)} \Biggl[\sum_{i_1 = 1}^\infty \la_{i_1}^{r \cdot 1} \Biggr]^{r_1} \Biggl[\sum_{i_2 = 1}^\infty \la_{i_2}^{r \cdot 2} \Biggr]^{r_2} \cdots \Biggl[\sum_{i_p = 1}^\infty \la_{i_p}^{r \cdot p} \Biggr]^{r_p} \\
& \quad = 2^p \Bigl[\tfrac12 \zeta_D (-r\cdot 1) \Bigr]^{r_1 } \Bigl[\tfrac12 \zeta_D (-r\cdot 2) \Bigr]^{r_2} \cdots \Bigl[\tfrac12 \zeta_D (-r\cdot p) \Bigr]^{r_p} .
\end{align*}
Therefore, with the sums below taken over nonnegative integers $r_1,\dots,r_p$ such that $1 r_1 + 2 r_2 + \cdots + p r_p = p$,
\begin{align*}
\mu_{r,F}(\si^2) &= \sigma^{2p} 2^p \sum_{} \Biggl[\begin{array}{c} p \\ r_1,\dots,r_p \end{array} \Biggr] [(1-1)!]^{r_1} [(2-1)!]^{r_2} \cdots [(p-1)!]^{r_p} \prod_{\ell=1}^p [\tfrac12 \zeta_D(-r\ell)]^{r_\ell} \\
&= \sigma^{2p} 2^p \sum_{} \Biggl[\begin{array}{c} p \\ r_1,\dots,r_p \end{array} \Biggr] \prod_{\ell=1}^p \Biggl[ \frac{(\ell-1)!}{2} \zeta_D (-r\ell) \Biggr]^{r_\ell} \\
&= \sigma^{2p} 2^p B_p(x_1,\dots,x_p),
\end{align*}
where $B_p$ is the complete Bell polynomial in $p$ variables  and $x_\ell$ is defined in \eqref{xl}.
\epr

\subsection{Lemmas used in the proof of Theorem \ref{LLN}, and their proofs}\label{nsec3.2}

   We now prove the lemmas that were quoted in the proof of Theorem \ref{LLN}.

\blem\label{remove-drift}
Let $r<-\frac d2$ and the assumptions of Theorem~\ref{LLN} hold. Let $u'(t,\cdot)$ be defined in \eqref{uprime}, where $X_k$ is given by \eqref{Xkt}. Then
\[ V^{n,r}_F(u,t)-V^{n,r}_F(u',t) \stackrel{L^1}{\longrightarrow} 0\qquad\text{as}\qquad n\to\infty. \]
\elem

For the proof of Lemma \ref{remove-drift}, we need an auxiliary result.
\blem\label{tight}
The family
\beq\label{calu} \calu :=\Biggl\{ \frac{\Delta^n_i u}{\sqrt{\Del}}, \frac{\Delta^n_i u'}{\sqrt{\Del}}: n\in\N,~ i=1,\ldots,[T/\Del]\Biggr\}  \eeq
of random elements in $H_r$ is tight.
\elem
\bpr We use the tightness criterion given in \cite[Theorem~1]{Suquet99} and consider the orthonormal basis in $H_r$ given by \eqref{basis_b}.
Then we have for any $1\leq m_1\leq m_2\leq \infty$,
\begin{align*}
\sum_{k=m_1}^{m_2} \langle \Delta^n_i u, b_k\rangle_{H_r}^2 &= \sum_{k=m_1}^{m_2} \Biggl( \sum_{\ell=1}^\infty \la_\ell^r  a_\ell(\Delta^n_i u)a_\ell(b_k)\Biggr)^2 = \sum_{k=m_1}^{m_2}\la_k^r \left[a_k(\Delta^n_i u)\right]^2 \\
&= \sum_{k=m_1}^{m_2}\la_k^r  \Biggl( \iint_0^{i\Del} \Delta^n_i \ee_k(s)\phi_k(y)\si(s,y)\,W(\dd s,\dd y)\Biggr)^2,
\end{align*}
and therefore by \eqref{mom} and calculations similar to those in \eqref{calc},
\beq\label{series-calc}\begin{split}
\bbe\Biggl[ \sum_{k=m_1}^{m_2} \langle \Delta^n_i u, b_k\rangle_{H_r}^2\Biggr] &= \sum_{k=m_1}^{m_2}\la_k^r \iint_0^{i\Del} (\Delta^n_i \ee_k(s))^2\phi_k(y)^2\bbe[\si^2(s,y)]\,\dd s\,\dd y\\
&\lec \sum_{k=m_1}^{m_2}\la_k^r \int_0^{i\Del} (\Delta^n_i \ee_k(s))^2\,\dd s\\
&=\sum_{k=m_1}^{m_2} \la_k^{r-\ga} (1-\ee^{-\la_k^\ga\Del}) - \frac12\sum_{k=m_1}^{m_2} \la_k^{r -\ga} \ee^{-2\la_k^\ga i\Del}  (\ee^{\la_k^\ga\Del}-1)^2 \\
&\leq \Del\sum_{k=m_1}^{m_2} \la_k^r. 
\end{split}\eeq
Hence, we obtain
\beq\label{m1m2}
\bbe\Biggl[ \sum_{k=m_1}^{m_2} \biggl\langle \frac{\Delta^n_i u}{\sqrt{\Del}}, b_k\biggr\rangle_{H_r}^2\Biggr] \lec \sum_{k=m_1}^{m_2} \la_k^r.
\eeq
The same estimate holds true when $\Delta^n_i u$ is replaced by $\Delta^n_i u'$ as the reader may quickly verify. We obtain on the one hand, for any $m\geq1$,
\[ \lim_{M\to\infty} \sup_{h \in \calu} \bbp\Biggl( \sum_{k=1}^m \langle h, b_k\rangle_{H_r}^2 > M \Biggr) \lec \lim_{M\to\infty}\frac1M \sum_{k=1}^{m} \la_k^r = 0, \]
and on the other hand, as a consequence of Lemma~\ref{A1} (i), for any $\delta>0$,
 \[ \lim_{m\to\infty} \sup_{h \in \calu} \bbp\Biggl( \sum_{k=m}^\infty \langle h, b_k\rangle_{H_r}^2 > \delta \Biggr) \lec \lim_{m\to\infty}\frac1\delta \sum_{k=m}^{\infty} \la_k^r = 0. \]
 Hence, the claim follows from \cite[Theorem~1]{Suquet99}.
\epr

\bpr[Proof of Lemma~\ref{remove-drift}]
By Lemma~\ref{tight}, we can find for every $\theta>0$, a compact subset $K_\theta$ of $H_r$ such that
\beq\label{theta}
\sup_{n\in\N} \sup_{i=1,\ldots,[t/\Del]} \Biggl\{ \bbp\biggl( \frac{\Delta^n_i u}{\sqrt{\Del}} \notin K_\theta\biggr) +\bbp\biggl( \frac{\Delta^n_i u'}{\sqrt{\Del}} \notin K_\theta\biggr) \Biggr\} <\theta.
\eeq
Moreover, for all $\theta>0$, since $F$ is uniformly continuous on $K_\theta$ by the Heine--Cantor theorem, there exist for given $\delta>0$, numbers $\eps_\theta(\delta)>0$ with the property $\eps_\theta(\delta)\to0$ as $\delta\to0$, such that
\beq\label{unifcont} h,h'\in K_\theta,~\|h-h'\|_{H_r} \leq \delta \implies |F(h)-F(h')|\leq \eps_\theta(\delta). \eeq
If we further let $M(\theta)=\sup_{h\in K_\theta} |F(h)|$, which is finite for all $\theta>0$, and recall that $|F(h)|\lec 1+|h|^p$ for some $p\geq2$, we derive
\beq\label{Fdiff} |F(h)-F(h')| \lec \eps_\theta(\delta)+M(\theta)\frac{\|h-h'\|_{H_r}}{\delta} + \bone_{\{h \text{ or } h'\notin K_\theta\}}(1+\|h\|_{H_r}^p + \|h'\|_{H_r}^p). \eeq

Suppose now that we can show the following:
\begin{align}
\label{toshow-1} \lim_{n\to\infty} \sup_{i=1,\ldots,[t/\Del]} \frac{\bbe[\|\Delta^n_i u-\Delta^n_i u'\|_{H_r}]}{\sqrt{\Del}} &= 0,\\
\label{toshow-2} \limsup_{n\to\infty} \sup_{i=1,\ldots,[t/\Del]} \frac{(\bbe[\|\Delta^n_i u\|_{H_r}^{2p}])^{1/2}+(\bbe[\|\Delta^n_i u'\|_{H_r}^{2p}])^{1/2}}{\Del^{p/2}}&<\infty.
\end{align}
Then using \eqref{Fdiff} and applying the Cauchy--Schwarz inequality to the third term results in
\begin{align*}
&\bbe[|V^{n,r}_F(u,t)-V^{n,r}_F(u',t)|]\\
&\quad\leq \Del\, [t/\Del]\sup_{i=1,\ldots,[t/\Del]} \bbe\Biggl[\biggl\lvert F\biggl(\frac{\Delta^n_i u}{\sqrt{\Del}}\biggr) -F\biggl( \frac{\Delta^n_i u'}{\sqrt{\Del}}\biggr) \biggr\rvert\Biggr]\\
&\quad\lec \eps_{\theta}(\delta) +M(\theta)\sup_{i=1,\ldots,[t/\Del]}\frac{\bbe[\|\Delta^n_i u-\Delta^n_i u'\|_{H_r}]}{\delta\sqrt{\Del}}+ \bbp\biggl(\frac{\Delta^n_i u}{\sqrt{\Del}}\text{ or }\frac{\Delta^n_i u'}{\sqrt{\Del}}\notin K_\theta\biggr)^{1/2}\\
&\qquad\quad\times\sup_{i=1,\ldots,[t/\Del]} \Biggl(1+\frac{(\bbe[\|\Delta^n_i u\|_{H_r}^{2p}])^{1/2}}{\Del^{p/2}} + \frac{(\bbe[\|\Delta^n_i u'\|_{H_r}^{2p}])^{1/2}}{\Del^{p/2}}\Biggr).
\end{align*}
Hence, \eqref{theta}, \eqref{toshow-1}, and \eqref{toshow-2} imply
\begin{align*}
\limsup_{n\to\infty} \bbe[|V^{n,r}_F(u,t)-V^{n,r}_F(u',t)|]\lec \eps_{\theta}(\delta)+\theta^{\frac12}.
\end{align*}
Letting $\delta\to0$ and then $\theta\to0$ yields the assertion of the lemma.

Therefore, it remains to verify \eqref{toshow-1} and \eqref{toshow-2}. For \eqref{toshow-2}, 
by the Minkowski and the Burkholder--Davis--Gundy inequalities and  calculations that are similar to \eqref{series-calc}, we have
\begin{align}\nonumber
(\E[\|\Delta^n_i u\|_{H_r}^{2p}])^{\frac1p} &= \Biggl(\bbe\Biggl[ \Biggl( \sum_{k=1}^\infty \la_k^r(\Delta^n_i a_k)^2\Biggr)^p\Biggr]\Biggr)^{1/p}\leq \sum_{k=1}^\infty \la_k^r\Bigl(\bbe[|\Delta^n_i a_k|^{2p}]\Bigr)^{\frac1p}\\
&\leq \sum_{k=1}^\infty \la_k^r \iint_0^{i\Del} (\Delta^n_i \ee_k(s))^2\phi_k(y)^2(\bbe[|\si(s,y)|^{2p}])^{\frac1p}\,\dd s\,\dd y\lec \Del,
\label{rd3.22}
\end{align}
from which the desired bound for $\Delta^n_i u$ follows. The arguments for $\Delta^n_i u'$ are completely analogous. This proves \eqref{toshow-2}.

For \eqref{toshow-1}, we first compute the variance of an increment of the drift process in \eqref{SDE}. By standard stochastic calculus and straightforward manipulations, we have
\begin{align*}
&\bbe[|\Delta^n_i a_k-\Delta^n_i X_k|^2]\\
&\quad= 2\la_k^{2\ga} \int_{(i-1)\Del}^{i\Del}\int_{(i-1)\Del}^s\bbe[a_k(s)a_k(r)]\,\dd r\,\dd s\\
&\quad=2 \la_k^{2\ga} \int_{(i-1)\Del}^{i\Del}\int_{(i-1)\Del}^s\iint_0^r \ee^{-\la_k^\ga(s-v)-\la_k^\ga(r-v)}\phi_k(y)^2\bbe[\si^2(v,y)]\,\dd v\,\dd y\,\dd r\,\dd s\\
&\quad\lec  2\la_k^{2\ga} \int_{(i-1)\Del}^{i\Del} \ee^{-\la_k^\ga s} \int_{(i-1)\Del}^s\ee^{-\la_k^\ga r}\int_0^r \ee^{2\la_k^\ga v}\,\dd v\,\dd r\,\dd s\\
& \quad=\la_k^\ga  \int_{(i-1)\Del}^{i\Del} \ee^{-\la_k^\ga s} \int_{(i-1)\Del}^s(\ee^{\la_k^\ga r}-\ee^{-\la_k^\ga r})\,\dd r\,\dd s\\
&\quad =\la_k^\ga  \int_{(i-1)\Del}^{i\Del} \ee^{-\la_k^\ga (s-(i-1)\Del)} \int_{(i-1)\Del}^s\ee^{-\la_k^\ga (i-1)\Del)}(\ee^{\la_k^\ga  r}-\ee^{-\la_k^\ga r})\,\dd r\,\dd s\\
&\quad =\la_k^\ga  \int_{(i-1)\Del}^{i\Del} \ee^{-\la_k^\ga (s-(i-1)\Del)} \int_0^{s-(i-1)\Del} (\ee^{\la_k^\ga r}-\ee^{-2\la_k^\ga (i-1)\Del)}\ee^{-\la_k^\ga r})\,\dd r\,\dd s\\
&\quad=\la_k^\ga  \int_0^{\Del} \ee^{-\la_k^\ga s} \int_0^{s} (\ee^{\la_k^\ga r}-\ee^{-2\la_k^\ga (i-1)\Del)}\ee^{-\la_k^\ga r})\,\dd r\,\dd s\\
&\quad=\int_0^{\Del} (1-\ee^{-\la_k^\ga s})\,\dd s - \ee^{-2\la_k^\ga (i-1)\Del}\int_0^{\Del} \ee^{-\la_k^\ga s}(1-\ee^{-\la_k^\ga s})\,\dd s\\
&\quad\leq \int_0^{\Del}(1-\ee^{-\la_k^\ga s})\,\dd s.
\end{align*}
Now let us choose a small $\rho\in(0,1)$ such that $r+\ga\rho<-\frac d2$. Then, with the estimate $1-\ee^{-\la_k^\ga  s}\leq (\la_k^\ga  s)^{\rho}(1-\ee^{-\la_k^\ga  s})^{1-\rho} \leq (\la_k^\ga  \Del)^{\rho}$ for all $s\leq\Del$, we obtain $\bbe[|\Delta^n_i a_k-\Delta^n_i X_k|^2] \lec \Del^{1+\rho}\la_k^{\ga \rho}$, and therefore,
\begin{align*}
\bbe[\|\Delta^n_i u - \Delta^n_i u'\|^2_{H_r}] &= \sum_{k=1}^\infty \la_k^r  \bbe[|\Delta^n_i a_k-\Delta^n_i X_k|^2]\lec \Del^{1+\rho}\sum_{k=1}^\infty \la_k^{r+\ga\rho} .
\end{align*}
By Lemma~\ref{A1} (i), the series converges with our choice of $\rho$, which proves \eqref{toshow-1}.
\epr

Due to the martingale property of $X_k$, the increments of $u'$ are uncorrelated to each other. However, we would like them to be more, namely, conditionally independent given each other. For this purpose, we have to discretize the random field $\si$.
\blem\label{discr-sigma}
Let $h^n_i$ be as defined in \eqref{hni}. As $n\to\infty$,
\[ V^{n,r}_F(u',t) - \Del\sumt F\biggl( \frac{h^n_i}{\sqrt{\Del}}\biggr) \stackrel{L^1}{\longrightarrow} 0. \]
\elem

\bpr
The proof follows the same scheme as the proof of Lemma~\ref{remove-drift}. The main tool is the estimate \eqref{Fdiff}, this time applied to $\Delta^n_i u'$ and $h^{n}_i$ (note that $\calu$ in \eqref{calu} remains tight if we include also the vectors $h^n_i/\sqrt{\Del}$). Therefore, the lemma is proved once we can show an analogue of \eqref{toshow-1} and \eqref{toshow-2}. We shall only give the details for the modified version of \eqref{toshow-1}.
Notice that, with the notation of \eqref{L2-cont},
\begin{align*}
\bbe[\|\Delta^n_i u'-h^n_i\|_{H_r}^2]&= \sum_{k=1}^\infty \la_k^r  \iint_{(i-1)\Del}^{i\Del} \phi_k(y)^2\bbe[|\si(s,y)-\si((i-1)\Del,y)|^2]\,\dd s\,\dd y\\
&\leq \Del\, w_2(\Del;t)^2 \sum_{k=1}^\infty \la_k^r.
\end{align*}
Since $w_2(\Del;t)\to0$ as $n\to\infty$, the statement corresponding to \eqref{toshow-1} follows.
\epr

The law of large numbers for $h^n_i$ can now be shown via covariance analysis.
\blem\label{LLN-core-1} We have
\beq\label{LLN-core-eq-1}   \lim_{n\to\infty} \bbe\Biggl[\Biggl\lvert \Del\sumt \Biggl\{  F\biggl( \frac{h^n_i}{\sqrt{\Del}}\biggr)-\bbe\biggl[ F\biggl( \frac{h^n_i}{\sqrt{\Del}}\biggr)\mathrel{\bigg|} \calf_{(i-1)\Del}\biggr]\Biggr\}\Biggr\rvert^2\Biggr] = 0. \eeq
\elem
\bpr By construction, if $j>i$, then $h^{n}_i$ and $h^{n}_j$ are conditionally independent given $\calf_{(j-1)\Delta_n}$. Moreover, since $F$ has at most polynomial growth, the variance of $F(h^n_i/\sqrt{\Del})$ can be bounded independently of $n$ and $i$, cf.\ \eqref{toshow-2}. Therefore, the second moment on the left-hand side of \eqref{LLN-core-eq-1} is bounded by a constant times $\Del^2\, [t/\Del]\leq t\Del$, which converges to $0$ as $n\to\infty$.
\epr

The next step is to compute the conditional expectation in \eqref{LLN-core-eq-1}. To this end, we need to recall and prove some results concerning Gaussian measures on Hilbert spaces. By \cite[Theorem~2.1.2]{Prevot07}, the second moment of a mean-zero Gaussian random vector on $H_r$ with covariance operator $Q$  is given by the \emph{trace-class norm} of $Q$:
\beq\label{trace} \|Q\|_1 := \sum_{k=1}^\infty \langle Qb_k,b_k\rangle_{H_r},  \eeq
where $(b_k)_{k\geq1}$ is an arbitrary orthonormal basis of $H_r$. Another norm on $\call^{+,\mathrm{sym}}_1(H_r)$ is given by the \emph{Hilbert--Schmidt norm}
\beq\label{HS} \|Q\|_2 := \Biggl(\sum_{k=1}^\infty \|Qb_k\|_{H_r}^2\Biggr)^{1/2}.\eeq
By \cite[Remark~B.0.4 and B.0.6]{Prevot07}, neither $\|Q\|_1$ nor $\|Q\|_2$ depends on the choice of $(b_k)_{k\geq1}$. We introduce a third metric on $\call^{+,\mathrm{sym}}_1(H_r)$ by setting
\beq\label{dr} d_r(Q,Q') := \|Q-Q'\|_{\call(H_r)} + | \|Q\|_1 - \|Q'\|_1|,\qquad Q,Q'\in \call^{+,\mathrm{sym}}_1(H_r),\eeq
where $\|\cdot\|_{\call(H_r)}$ is the usual operator norm for bounded linear operators $H_r\to H_r$.

\blem\label{comp-condexp}
We have
\[ \bbe\Biggl[ F\biggl( \frac{h^n_i}{\sqrt{\Del}}\biggr)\mathrel{\bigg|} \calf_{(i-1)\Del}\Biggr] = \mu_F(Q^{n,i}_r), \]
where $Q^{n,i}_r\colon H_r \to H_r$ is given by
\beq\label{Qni}
(Q^{n,i}_r h)(h'):= \sum_{k,\ell=1}^\infty  \la_\ell^r  a_k(h')a_\ell(h)\int_D \phi_k(y)\phi_\ell(y) \si^2((i-1)\Del,y)\,\dd y\eeq
for $h\in H_r$ and $h'\in H_{-r}$.
\elem
\bpr For any $h\in H_r$, we have
\[
\langle h^n_i, h\rangle_{H_r} = \sum_{k=1}^\infty \la_k^r  a_k(h^n_i)a_k(h) = \sum_{k=1}^\infty \la_k^r a_k(h) \iint_{(i-1)\Del}^{i\Del}\phi_k(y)\si((i-1)\Del,y)\,W(\dd s,\dd y),
\]
which is, conditionally on $\calf_{(i-1)\Del}$, normally distributed with mean $0$ and variance
\beq\label{vni} v^n_i(h):= \Del\sum_{k,\ell=1}^\infty \la_k^r\la_\ell^r a_k(h)a_\ell(h)c^{n,i}_{k,\ell}, \eeq
where $c^{n,i}_{k,\ell}$ is the integral in \eqref{Qni}. Hence, by \cite[Definition~2.1.1]{Prevot07}, $h^n_i/\sqrt{\Del}$ is conditionally Gaussian given $\calf_{(i-1)\Del}$ with zero mean. It remains to show that its conditional covariance operator is given by $Q^{n,i}_r$. But this follows
immediately from \eqref{Qni} because
\beq\label{ak}\begin{aligned}
	a_k(Q^{n,i}_r h) = (Q^{n,i}_r h)(\phi_k)= \sum_{\ell=1}^\infty\la_\ell^r a_\ell(h) c^{n,i}_{k,\ell}
\end{aligned}\eeq
and consequently,
\beq\label{QHH} \langle Q^{n,i}_r h, h\rangle_{H_r} = \sum_{k=1}^\infty \la_k^r a_k(Q^{n,i}_r h)a_k(h)=\frac{1}{\Del} v^n_i(h). \eeq
\epr

\blem\label{muf-prop}
If $F\colon H_r\to\R$ is continuous with $|F(h)|\lec 1+\|h\|_{H_r}^p$ for some $p>0$, then $\mu_F$ is continuous on $\call^{+,\mathrm{sym}}_1(H_r)$ with respect to $d_r$, and $|\mu_F(Q)| \lec 1+\|Q\|_{1}^{p/2}$.
\elem
\bpr Since $\|Q\|_1$ is precisely the second moment of $N_r(0,Q)$, the bound on $\mu_F(Q)$ follows immediately from the growth assumption on $F$ and \cite[Corollary~3.2]{Ledoux91}. For the continuity of $\mu_F$, let
$Q_n$ and $Q$ be elements of $\call^{+,\mathrm{sym}}_1(H_r)$ such that $d_r(Q_n,Q)\to0$.
Then, the assumptions of \cite[Example~3.8.15]{Bogachev98} are satisfied, and we have $N_r(0,Q_n)\limw N_r(0,Q)$ as $n\to\infty$, where $\limw$ denotes weak convergence of measures.
Since $F$ is of polynomial growth and the moments of any order of $N_r(0,Q_n)$ are uniformly bounded in $n$ (as the first part of the proof has shown), we obtain $\mu_F(Q_n)\to \mu_F(Q)$ from \cite[Lemma 3.8.7]{Bogachev98}.
\epr

The next lemma is the final ingredient in the proof of Theorem~\ref{LLN} for $r<-\frac d2$.
\blem\label{cex-conv} We have
\beq\label{eq1} \Del\sumt \mu_F(Q^{n,i}_r) - \int_0^t \mu_F(\si^2(s,\cdot))\,\dd s \stackrel{L^1}{\longrightarrow}0\qquad\text{as} \qquad n\to\infty. \eeq
\elem
\bpr Let $Q_r(v)$ (resp., $c_{k,\ell}(v)$) be defined in the same way as $Q^{n,i}_r$ in \eqref{Qni} (resp., the integral in \eqref{Qni}), but with $(i-1)\Del$ replaced by $v$. By definition (see \eqref{Qphi1}), we have $\mu_F(Q_r(v))=\mu_{r,F}(\si^2(v,\cdot))$ and $Q^{n,i}_r=Q_r((i-1)\Del)$. Therefore, we can write the difference in \eqref{eq1} as $C^n_1+C^n_2$ where
\begin{align*}
C^n_1 &:=\int_0^t \Biggl\{\sumt \Bigl(\mu_F(Q_r((i-1)\Del)) - \mu_F(Q_r(v))\Bigr)\bone_{((i-1)\Del,i\Del]}(v)\Biggr\}\, \dd v,\\
C^n_2&:=- \int_{[t/\Del]\Del}^{t} \mu_F(Q_r(v))\,\dd v.
\end{align*}
Suppose we show that
\beq\label{Q-conv} \sup_{v,w\in[0,t]:\,|v-w|\leq\Del}\bbe[d_r(Q_r(v),Q_r(w))^2] \lec w_2(\Del;t)^2. \eeq
Then, together with the $L^2(\Om)$-continuity of $\si$ and the continuity of $Q\mapsto \mu_F(Q)$ (see Lemma~\ref{muf-prop}), this will imply that for every $v\in[0,t]$, the integrand in the definition of $C^n_1$ converges to $0$ in probability as $n\to\infty$. In fact, this convergence also takes place in $L^1(\Om)$ by uniform integrability: indeed, by Lemma~\ref{muf-prop}, \eqref{mom} and equations \eqref{trace}, \eqref{vni}, and \eqref{QHH},
\beq\label{muf-2}\begin{aligned}
\bbe[\mu_F(Q_r(v))^2] &\lec 1+ \bbe[\|Q_r(v)\|_1^p] \leq 1+\bbe\Biggl[ \Biggl(  \sum_{k=1}^\infty \la_k^rc_{k,k}(v)\Biggr)^p\Biggr]\\
&\leq 1+\Biggl(  \sum_{k=1}^\infty \la_k^r \bbe[|c_{k,k}(v)|^p]^{\frac1p}\Biggr)^p \lec 1+\Biggl(\sum_{k=1}^\infty \la_k^r\Biggr)^p,
\end{aligned}\eeq
which is finite and independent of $v$. Therefore, for $v \in [0,t]$,
\begin{align*}
    &  \bbe\Biggl[ \Biggl(\sumt \Bigl(\mu_F(Q_r((i-1)\Del) - \mu_F(Q_r(v))\Bigr)\bone_{((i-1)\Del,i\Del]}(v)\Biggr)^2 \Biggr] \\
   &\quad =  \bbe \Biggl[\sumt \Bigl(\mu_F(Q_r((i-1)\Del) - \mu_F(Q_r(v))\Bigr)^2\bone_{((i-1)\Del,i\Del]}(v)\Biggr]  \\
   &\quad\lec 2 \sumt \Biggl( 1+\Biggl(\sum_{k=1}^\infty \la_k^r\Biggr)^p \Biggr) \bone_{((i-1)\Del,i\Del]}(v)   \leq 2 \Biggl( 1+\Biggl(\sum_{k=1}^\infty \la_k^r\Biggr)^p\Biggr)  < \infty.
\end{align*}
This establishes the uniform integrability of the integrand in the definition of $C^n_1$. Applying dominated convergence to the $\dd v$-integral, we obtain $\bbe[|C^n_1|]\to0$. The same estimate \eqref{muf-2} also shows $\bbe[|C^n_2|]\to0$.

Hence, it remains to prove \eqref{Q-conv}, and we shall consider the two parts defining the metric $d_r$ in \eqref{dr} separately. Using \cite[Remark~B.0.6(ii)]{Prevot07}, we obtain
\beq\label{eq2}
\|Q_r(v)-Q_r(w)\|_{\call(H_r)}^2 \leq \|Q_r(v)-Q_r(w)\|_2^2 = \sum_{k=1}^\infty \| Q_r(v)b_k-Q_r(w)b_k\|_{H_r}^2.
\eeq
Recalling the formula for $b_k$ from \eqref{basis_b},
\[
\langle Q_r(v)b_k,Q_r(w)b_k\rangle_{H_r} = \sum_{j=1}^\infty \la_j^r a_j(Q_r(v)b_k)a_j(Q_r(w)b_k)= \la_k^r \sum_{j=1}^\infty \la_j^r c_{j,k}(v)c_{j,k}(w)
\]
by \eqref{ak}, so expanding the square in the last term in \eqref{eq2} yields for $v,w \in [0,1]$ with $\lvert v - w \rvert \leq \Del$,
\begin{align*}
\bbe[\|Q_r(v)-Q_r(w)\|_{\call(H_r)}^2]\leq   \sum_{k=1}^\infty \la_k^r \sum_{j=1}^\infty \la_j^r \bbe[ (c_{j,k}(v)-c_{j,k}(w))^2] \lec w_2(\Del;t)^2
\end{align*}
where the last inequality uses the Cauchy-Schwartz inequality and the fact that $\sigma$ has bounded moments by \eqref{mom}, as well as Lemma~\ref{A1} (i)  and the fact that $r < -\tfrac d2$. For the second part of the metric $d_r$, observe from \eqref{vni} and \eqref{QHH} that
\begin{align*}
\|Q_r(v)\|_1-\|Q_r(w)\|_1 &= \sum_{k=1}^\infty \Big( \langle Q_r(v)b_k,b_k\rangle_{H_r} -  \langle Q_r(w)b_k,b_k\rangle_{H_r}\Big)\\
&= \sum_{k=1}^\infty \la_k^r \Big( c_{k,k}(v) - c_{k,k}(w)\Big).
\end{align*}
This in turn gives
\begin{align*}
\bbe\biggl[\Big|\|Q_r(v)\|_1-\|Q_r(w)\|_1\Big|^2\biggr] \leq \Biggl(  \sum_{k=1}^\infty \la_k^r \bbe[( c_{k,k}(v) - c_{k,k}(w))^2]^{\frac12} \Biggr)^2 \lec w_2(\Del;t)^2,
\end{align*}
which completes the proof of \eqref{Q-conv}.
\epr


\subsection{Proof of Theorem~\ref{LLN2} ($-\tfrac{d}{2} \leq r <\ga- \tfrac d2$)}\label{nsec3.3}
\bpr[Proof of Theorem~\ref{LLN2}]
As in the proof of Theorem~\ref{LLN}, we will approximate $\|\Delta^n_i u\|_{H_r}$ by simpler expressions: For $\eps>0$, $i= 1,\dots, [t/\Del]$, $k \geq 1$, and $-\frac d2\leq r<\ga-\frac d2$, we consider, in a first truncation step,
\begin{equation}\label{hnieps} a^{n,i,\eps}_k:=\iint_{(i\Del-\eps)\vee 0}^{i\Del} \Delta^n_i \ee_k(s)\phi_k(y)\si(s,y)
\,W(\dd s,\dd y),\qquad
h^{n,i,\eps}_r :=\Biggl(\sum_{k=1}^\infty \la_k^r(a^{n,i,\eps}_k)^2\Biggr)^{1/2}\end{equation}
as approximations for $\Delta^n_i a_k$ (as defined in \eqref{abbr}) and $\|\Delta^n_i u\|_{H_r}$, respectively. In a second step, we discretize the random field $\si$ by introducing the variables
\beq\label{wthnieps}\begin{split} \wt a^{n,i,\eps}_k&:=\iint_{(i\Del-\eps)\vee 0}^{i\Del} \Delta^n_i \ee_k(s)\phi_k(y)\si(i\Del-\eps,y)
\,W(\dd s,\dd y),\\
\wt h^{n,i,\eps}_r  &:=\Biggl(\sum_{k=1}^\infty \la_k^r(\wt a^{n,i,\eps}_k)^2\Biggr)^{1/2},  \end{split}\eeq
where  $\si$ is extended to $\R\times D$ by setting $\si(s,y):=\si(0,y)$ for all $s<0$ and $y\in D$.

In this set-up, we consider the decomposition
\[ V^{n,r}_f(u,t)-V^r_f(u,t)=J^\eps_{n,1}+J^\eps_{n,2}+J^\eps_{n,3},  \]
where
\begin{align*}
J^\eps_{n,1}&:=V^{n,r}_f(u,t)-\Del\sumt f\biggl( \frac{h^{n,i,\eps}_r}{\tau_n(r)}\biggr), \\
J^\eps_{n,2}&:=\Del\sumt   f\biggl( \frac{h^{n,i,\eps}_r}{\tau_n(r)}\biggr)-\Del\sumt f\biggl( \frac{\wt h^{n,i,\eps}_r}{\tau_n(r)}\biggr),\\
J^\eps_{n,3}&:=\Del\sumt  f\biggl( \frac{\wt h^{n,i,\eps}_r}{\tau_n(r)}\biggr) -V^r_f(u,t)
\end{align*}
and prove \eqref{eq:LLN-2} by showing
\[ \lim_{\eps\to0}\limsup_{n\to\infty} \E[|J^\eps_{n,1}|+|J^\eps_{n,2}|+|J^\eps_{n,3}|] =0 \]
in Lemmas~\ref{lemma-eps}, \ref{LLN-discr}, and \ref{LLN-condexp} below. Indeed, as mentioned at the beginning of Section \ref{sec3.1}, having $\E[|V^{n,r}_f(u,t)-V^r_f(u,t)|]\to0$ for all $t\in[0,T]$  is sufficient to obtain the uniform convergence to $0$ on $[0,T]$ in $L^q(\Om)$ that is claimed in \eqref{eq:LLN-2}. This establishes Theorem \ref{LLN2}.
\epr

While Lemmas~\ref{lemma-eps} and \ref{LLN-discr} (needed for the proof of Theorem~\ref{LLN}) have similar proofs to Lemmas~\ref{remove-drift} and \ref{discr-sigma} (needed for the proof of Theorem~\ref{LLN2}), both the assertion and the proof of Lemma~\ref{LLN-condexp} are significantly different from what we saw in Lemma~\ref{LLN-core-1}. In fact, the main step in proving Lemma~\ref{LLN-condexp} is Lemma~\ref{convP}, where we show that, in the case $-\frac d2\leq r<\ga-\frac d2$, not only the  \emph{average} but actually \emph{each} of the variables $\wt h^{n,i,\eps}_r/\tau_n(r)$ converges in $L^2(\Om)$; cf.\ Remark~\ref{difference} (ii).

\subsection{Lemmas used in the proof of Theorem \ref{LLN2}, and their proofs}\label{nsec3.4}
\blem\label{lemma-eps} Let $-\frac d2\leq r<\ga-\frac d2$ and consider the variables $h^{n,i,\eps}_r$ from \eqref{hnieps}.
Then
\beq\label{diff} \lim_{\eps\to0} \limsup_{n\to\infty} \bbe\Biggl[\biggl\lvert V^{n,r}_f(u,t)-\Del\sumt f\biggl( \frac{h^{n,i,\eps}_r}{\tau_n(r)}\biggr)\biggr\rvert\Biggr] = 0. \eeq
\elem
\bpr The difference in \eqref{diff} is
\begin{align*}
D^{n,\eps}_r &:=\Del\sum_{i=1}^{[t/\Del]} \Biggl\{f\biggl(\frac{\|\Delta^n_i u\|_{H_r}}{\tau_n(r)}\biggr) -f\biggl( \frac{h^{n,i,\eps}_r}{\tau_n(r)}\biggr)\Biggr\}
\end{align*}
(note that for $i=1,\dots, [\eps/\Del]$, the terms in this sum vanish).
Since $f$ is continuous with polynomial growth, we can find $p\geq2$ such that $|f(x)|/x^p\to 0$ as $x\to \infty$. Hence, similarly to \cite[Equation~(3.4.16)]{Jacod12}, we have
\beq\label{fdiff} |f(x+y)-f(x)| \leq C\biggl( \Phi'_A(\delta) + \Phi(A) \frac{|y|}{\delta} + \Phi''(A)(x^p+y^p)\biggr),\qquad x,y\geq0,\eeq
where $C>0$ only depends on $p$, and $\Phi(A)$, $\Phi'_A(\delta)$, and $\Phi''(A)$ are finite numbers satisfying $\Phi'_A(\delta)\to0$ as $\delta\to0$ for given $A>0$ and $\Phi''(A)\to0$ as $A\to\infty$. If we can show that
\begin{align}
\label{toshow1} &\lim_{n\to\infty} \sup_{i=1,\ldots,[t/\Del]} \frac{\bbe[|\|\Delta^n_i u\|_{H_r}-h^{n,i,\eps}_r|]}{\tau_n(r)} = 0,\\
\label{toshow2} &\limsup_{\eps\to0} \limsup_{n\to\infty} \sup_{i=1,\ldots,[t/\Del]} \frac{\bbe[\|\Delta^n_i u\|_{H_r}^p]+\bbe[(h^{n,i,\eps}_r)^p]}{\tau_n(r)^p}<\infty,
\end{align}
then \eqref{fdiff} implies
\begin{align*}
\limsup_{n\to\infty} \bbe[|D^{n,\eps}_r|] &\leq \limsup_{n\to\infty} \Del\, [t/\Del]\sup_{i=1,\ldots,[t/\Del]} \bbe\Biggl[\biggl\lvert f\biggl(\frac{\|\Delta^n_i u\|_{H_r}}{\tau_n(r)}\biggr) -f\biggl( \frac{h^{n,i,\eps}_r}{\tau_n(r)}\biggr) \biggr\rvert\Biggr]\\
&\lec \Phi'_A(\delta)+\Phi''(A),
\end{align*}
which goes to $0$ if we first let $\delta\to0$ and then $A\to\infty$. 
This establishes \eqref{diff}.

So it remains to show \eqref{toshow1} and \eqref{toshow2}.
For \eqref{toshow2}, we note that
\begin{align*}
   \bbe[\|\Delta^n_i u\|_{H_r}^p] &= \E\Biggl[\Biggl(\sum_{k=1}^\infty \la_k^r (\Delta^n_i a_k)^2 \Biggr)^{p/2} \Biggr] = \Biggl\lVert \sum_{k=1}^\infty \la_k^r(\Delta^n_i a_k)^2 \Biggr\rVert_{L^{p/2}(\Omega)}^{p/2}  \\
   &\leq \Biggl(\sum_{k=1}^\infty \la_k^r \lVert \Delta^n_i a_k \rVert_{L^p(\Omega)}^2 \Biggr)^{p/2}
\end{align*}
by Minkowski's inequality. Now we apply the Burkholder--Davis--Gundy inequality, and again the Minkowski's inequality and hypothesis \eqref{mom} to see, as in the \eqref{calc}, that
\begin{align*}
   \lVert \Delta^n_i a_k \rVert_{L^p(\Omega)}^2 &\leq \Biggl\lVert \iint_0^{i\Del} (\Delta^n_i \ee_k(s))^2\phi_k(y)^2\si^{2}(s,y)\,\dd s\,\dd y \Biggr\rVert_{L^{p/2}(\Omega)} \\
    &\leq \iint_0^{i\Del} (\Delta^n_i \ee_k(s))^2\phi_k(y)^2(\bbe[|\si(s,y)|^p])^{\frac 2p}\,\dd s\,\dd y \lec \frac{1-\ee^{-\la_k^\ga\Del}}{\la_k^\ga},
\end{align*}
so
$
   \bbe[\|\Delta^n_i u\|_{H_r}^p] \lec  (A^{n,1}_r)^{p/2},
$
where $A^{n,1}_r/\tau_n(r)^2 \to K_r$ as seen in the proof of Lemma~\ref{size-incr}.
This show \eqref{toshow2} for the first term.  The proof for the second term is similar and is left to the reader.

We now prove \eqref{toshow1}. To this end, we use the reverse triangle inequality to deduce for $i\Del\geq\eps$,
\beq\label{diff-mom}\begin{aligned}
	\bbe[(\|\Delta^n_i u\|_{H_r}-h^{n,i,\eps}_r)^2]&\leq \sum_{k=1}^\infty \la_k^r \bbe[(\Delta^n_i a_k- a^{n,i,\eps}_k)^2]\\
	&=\sum_{k=1}^\infty \la_k^r \iint_0^{i\Del-\eps} (\Delta^n_i \ee_k(s))^2\phi_k(y)^2\bbe[\si^2(s,y)]\,\dd s\,\dd y\\
	&\lec \sum_{k=1}^\infty \la_k^r \int_0^{i\Del-\eps} (\ee^{-\la_k^\ga(i\Del-s)}-\ee^{-\la_k^\ga((i-1)\Del-s)})^2\,\dd s\\
	&=\frac12\sum_{k=1}^\infty \la_k^{r-\ga}(1-\ee^{\la_k^\ga\Del})^2(\ee^{-2\la_k^\ga\eps}-\ee^{-2\la_k^\ga i\Del})\\
	&\leq\frac12 \sum_{k=1}^\infty \la_k^{r-\ga}(1-\ee^{\la_k^\ga\Del})^2\ee^{-2\la_k^\ga\eps}.
\end{aligned}\eeq
Due to the presence of the exponential term $\ee^{-2\la_k^\ga\eps}$,
the right-hand side converges to a finite limit after dividing by $\Delta_n^2$, which shows \eqref{toshow1} by Jensen's inequality because $\Delta_n/\tau_n(r)\to0$.
\epr

\blem\label{LLN-discr}
With the variables $h^{n,i,\eps}_r$ and $\wt h^{n,i,\eps}_r$ from \eqref{hnieps} and \eqref{wthnieps}, we have,  as $n\to\infty$, that
\beq\label{diff2} \lim_{\eps\to0} \limsup_{n\to\infty} \Del\sumt \bbe\Biggl[\biggl\lvert  f\biggl( \frac{h^{n,i,\eps}_r}{\tau_n(r)}\biggr)-f\biggl( \frac{\wt h^{n,i,\eps}_r}{\tau_n(r)}\biggr)\biggr\rvert\Biggr] = 0. \eeq
\elem
\bpr
As in the proof of Lemma~\ref{lemma-eps}, it suffices to show \eqref{toshow1} and \eqref{toshow2} with $\|\Delta^n_i u\|_{H_r}$ replaced by $\wt h^{n,i,\eps}_r$. We shall only show the modified version of \eqref{toshow1}. With the same calculations as in \eqref{diff-mom}, we obtain for all $i\geq 1$,
\begin{align*}
&\bbe[(h^{n,i,\eps}_r-\wt h^{n,i,\eps}_r)^2] \\
&\quad\leq \sum_{k=1}^\infty  \la_k^r \iint_{(i\Del-\eps)\vee 0}^{i\Del} (\Delta^n_i \ee_k(s))^2\phi_k(y)^2\bbe[(\si(s,y)-\si(i\Del-\eps,y))^2]\,\dd s\,\dd y\\
&\quad\lec w_2(\eps;t)^2 \sum_{k=1}^\infty  \la_k^r \Biggl( \int_0^{i\Del} \ee^{-2 \la_k^\ga s} \,\dd s+ \int_0^{(i-1)\Del}\ee^{-2 \la_k^\ga s} \,\dd s\Biggr) \lec w_2(\eps;t)^2 \sum_{k=1}^\infty \la_k^{r-\ga},
\end{align*}
so the $L^2(\Om)$-continuity of $\si$ completes the proof since $r -\ga< -\frac d2$.
\epr


\blem\label{Pinr} For $n\in\N$ and $-\frac d2\leq r<\ga-\frac d2$, consider the measure on $[0,\infty)\times D$ defined by
\beq\label{Pin} \Pi^{n}_{r,\ga}(\dd s,\dd y) := \frac1{\tau_n(r)^2} \sum_{k=1}^\infty \la_k^r (\ee^{-\la_k^\ga s}-\ee^{-\la_k^\ga(s-\Del)}\bone_{s\geq\Del})^2\phi_k(y)^2\,\dd s\,\dd y. \eeq
Then, as $n\to\infty$,
\beq\label{weak} \Pi^{n}_{r,\ga} \limw \frac{K_r}{|D|}\,\delta_0(\dd s)\, \bone_D(y)\,\dd y, \eeq
where $K_r$ is the number from \eqref{Kr}.
\elem
\bpr Clearly, \eqref{weak} is equivalent to showing
\beq\label{weak-2} \begin{split} &\lim_{n\to\infty} \Pi^n_{r,\ga}([0,t]\times D_z) \\
	&\quad= \lim_{n\to\infty} \int_0^t \int_{D_z} \frac1{\tau_n(r)^2} \sum_{k=1}^\infty \la_k^r (\ee^{-\la_k^\ga s}-\ee^{-\la_k^\ga (s-\Del)}\bone_{s\geq\Del})^2\phi_k(y)^2\,\dd s\,\dd y \\
	&\quad= K_r\frac{|D_z|}{|D|} = \frac{K_r}{|D|}\iint \bone_{[0,t]\times D_z}(s,y)\,\delta_0(\dd s)\,\dd y\end{split}\eeq
for all $t> 0$ and $z\in \R^d$, where $D_z =\{y\in D: y_i\leq z_i \text{ for all } i=1,\dots,d\}$.

Consider the double integral in the first line of \eqref{weak-2}. With similar calculations as in \eqref{calc}, if we first integrate with respect to $s$, it is equal to
\beq\label{weak-3}  \frac1{\tau_n(r)^2}  \int_{D_z} \sum_{k=1}^\infty \la_k^{r-\ga}\Biggl( (1-\ee^{-\la_k^\ga\Del})-\frac{\ee^{-2\la_k^\ga t}}{2}(\ee^{\la_k^\ga \Del}-1)^2\Biggr)\phi_k(y)^2\,\dd y.  \eeq
Thus, \eqref{weak-2} is equivalent to the convergence of the previous line to $\frac{K_r}{|D|}\int_{D_z} 1\, \dd y$ for all $t>0$ and $z\in\R^d$. If $z$ is such that $D_z=D$, then \eqref{weak-3} is exactly the last line in \eqref{calc} (with $i\Del$ replaced by $t$) divided by $\tau_n(r)^2$, which, as seen in the proof of Lemma~\ref{size-incr}, converges to $K_r$. As a result, by Scheffé's theorem (see \cite[Theorem~16.12]{Billingsley95}), the lemma is proved if we can show that
\[ \frac1{\tau_n(r)^2}   \sum_{k=1}^\infty \la_k^{r-\ga}\Biggl( (1-\ee^{-\la_k^\ga\Del})-\frac{\ee^{-2\la_k^\ga t}}{2}(\ee^{\la_k^\ga \Del}-1)^2\Biggr)\phi_k(y)^2 \to \frac{K_r}{|D|} \]
as $n\to\infty$, pointwise for all $y\in D$. But this follows from Lemma~\ref{convdens}
because the second term inside the parentheses is negligible as the calculations for $A^{n,i,2}_r$ in \eqref{second-term} demonstrate.
\epr

The next lemma is the main step in the proof of Theorem~\ref{LLN2}. 
\blem\label{convP}
For $n\in\N$, $k\in\N$, $\eps>0$, $-\frac d2\leq r<\ga-\frac d2$, and $v\geq\eps$, define
\beq\label{tildeahu}\begin{aligned} \wt a^{n,\eps}_k(v)&:= \iint_{v-\eps}^{v} (\ee^{-\la_k^\ga(v-s)}-\ee^{-\la_k^\ga(v-\Del-s)}\bone_{s\leq v-\Del})\phi_k(y)\si(v-\eps,y)\,W(\dd s,\dd y),\\
	\wt h^{n,\eps}_r(v) &:=\Biggl(\sum_{k=1}^\infty \la_k^r(\wt a^{n,\eps}_k(v))^2\Biggr)^{1/2}. \end{aligned}\eeq
Then, for any $q>0$,
\beq\label{convP-eq} \sup_{v\in[\eps,t]}\bbe\Biggl[\Biggl\lvert \frac{\wt h^{n,\eps}_r(v)}{\tau_n(r)} - \biggl(\frac{K_r}{|D|}\int_D \si^2(v-\eps,y)\,\dd y\biggr)^{1/2}\Biggr\rvert^q\Biggr] \to0\as. \eeq
\elem
\bpr For simplicity, we adopt the notation $\bbe_{v-\eps}[\cdot] = \bbe[\cdot\,|\,\calf_{v-\eps}]$, and similarly $\var_{v-\eps}$ and $\cov_{v-\eps}$. As the moments of any order of $\wt h^{n,\eps}_r(v) /\tau_n(r)$ are uniformly bounded in $n$ and $v$ (cf.\ \eqref{toshow2}) and $\si^2(v-\eps,y)$ has uniformly bounded moments of all orders by \eqref{mom}, it suffices to prove the statement for $q=1$. Therefore, the result will follow from Lemma \ref{lemapp}  (with $\calg = \calf_{v-\eps}$) if we can show that, as $n\to\infty$,
\begin{align}
\label{conv-mean} \sup_{v\in[\eps,t]} \bbe\Biggl[\biggl\lvert\tau_n(r)^{-2}\bbe_{v-\eps}[(\wt h^{n,\eps}_r(v))^2] - \frac{K_r}{|D|}\int_D \si^2(v-\eps,y)\,\dd y\biggr\rvert\Biggr]&\to0,\\
\sup_{v\in[\eps,t]} \bbe\Bigl[\tau_n(r)^{-4}\var_{v-\eps}[(\wt h^{n,\eps}_r(v))^2]\Bigr]&\to0.
\label{conv-var}
\end{align}

Given a large number $L>0$, consider a partition $D=\bigcup_{\ell=1}^{L} D^L_\ell$ of $D$ into  pairwise disjoint measurable sets $D^L_l$ such that
\[ \lim_{L\to\infty} \sup_{\ell=1,\dots,L} \sup_{x,y\in D^L_\ell} |x-y| =0, \]
where $|\cdot|$ is an arbitrary norm on $\R^d$.
Picking an arbitrary point $y^L_\ell$ from each $D^L_\ell$, we can then
decompose the difference in \eqref{conv-mean}   into $E^{n,\eps,L,1}_r(v)+E^{n,\eps,L,2}_r(v)+E^{\eps,L,3}_r(v)$, where
\beq\label{D}\begin{aligned}
E^{n,\eps,L,1}_r(v)&:= \iint_0^v  \sum_{\ell=1}^L (\si^2(v-\eps,y)-\si^2(v-\eps,y^L_\ell))\bone_{D^L_\ell}(y)\bone_{s<\eps}\,\Pi^n_{r,\ga}(\dd s,\dd y),\\
E^{n,\eps,L,2}_r(v)&:=\sum_{\ell=1}^L \si^2(v-\eps,y^L_\ell) (\Pi^{n}_{r,\ga}([0,\eps]\times D^L_\ell)-\tfrac{K_r}{|D|}|D^L_\ell|),\\
E^{\eps,L,3}_r(v)&:=\frac{K_r}{|D|}\sum_{\ell=1}^L \int_D (\si^2(v-\eps,y^L_\ell) -\si^2(v-\eps,y))\bone_{D^L_\ell}(y)\,\dd y.
\end{aligned}\eeq
By \eqref{mom} and Lemma~\ref{Pinr}, we have $\sup_{v\in[\eps,t]}\bbe[|E^{n,\eps,L,2}_r(v)|]\to0$  as $n\to\infty$. Using the $L^2(\Om)$-continuity of $\si$, one can further show that
\[ \sup_{v\in[\eps,t]}\bbe[|E^{\eps,L,3}_r(v)|]\to 0\qquad\text{and}\qquad \limsup_{n\to\infty} \sup_{v\in[\eps,t]} \bbe[|E^{n,\eps,L,1}_r(v)|]\to 0\] as $L\to\infty$, proving \eqref{conv-mean}.

Next, we turn to \eqref{conv-var}. Using the notation $\Delta^n\ee_k(s)=\ee^{-\la_k^\ga s}-\ee^{-\la_k^\ga(s-\Del)}\bone_{s\geq\Del}$, and the formula $\cov(X^2,Y^2)=2\cov(X,Y)^2$ for mean-zero bivariate Gaussian variables $X$ and $Y$ (which follows, for example, from \eqref{prodmom}), we derive
\beq\label{var-line1}\begin{aligned}
\frac{1}{\tau_n(r)^4}\var_{v-\eps}[(\wt h^{n,\eps}_r(v))^2] &=\frac{1}{\tau_n(r)^4}\sum_{k,m=1}^\infty \la_k^r \la_m^r\cov_{v-\eps}\Bigl( \wt a^{n,\eps}_k(v)^2,\wt a^{n,\eps}_m(v)^2\Bigr)\\
&=\frac{2}{\tau_n(r)^4}\sum_{k,m=1}^\infty \la_k^r \la_m^r\biggl(\iint_0^\eps \Delta^n\ee_k(s)\Delta^n\ee_m(s) \\
&\qquad\times \phi_k(y)\phi_m(y)\si^2(v-\eps,y)\,\dd s\,\dd y\biggr)^2.
\end{aligned}\eeq

As $D$ is a bounded open set, we know
from the Vitali covering theorem (see \cite[Theorem~5.5.2]{Bogachev98}) that for small $\delta>0$, there are $L(\delta)\in\N$ and pairwise disjoint open balls $(B^\delta_\ell)_{\ell=1,\ldots, L(\delta)}$ of radius $\leq\delta$ such that $B^\delta_\ell \subseteq D$ and $\dist(B^\delta_\ell, \partial D)>\delta$ for every $i$ and such that $|D\setminus\bigcup_{\ell=1}^{L(\delta)} B^\delta_\ell| \downarrow 0$ as $\delta\downarrow 0$.
Denoting the center of $B^\delta_\ell$ by $z^\delta_\ell$ 
and using the elementary estimate $(x_1+x_2+x_3)^2\leq 3(x_1^2+x_2^2+x_3^2)$, we bound the terms in  \eqref{var-line1} by $3\sum_{i=1}^3 F^{n,\eps,\delta,i}_r(v)$, where
\begin{align*}
F^{n,\eps,\delta,1}_r(v)  &:=\frac{2}{\tau_n(r)^4}\sum_{k,m=1}^\infty \la_k^{r}\la_m^{r} \Biggl(\int_0^{\eps}  \Delta^n\ee_k(s)\Delta^n\ee_m(s) \,\dd s\\
&\qquad\times \int_D \phi_k(y)\phi_m(y)  \sum_{\ell=1}^{L(\delta)}\bone_{B^\delta_\ell}(y)  \Bigl(\si^2(v-\eps,y) -\si^2(v-\eps,z^\delta_\ell) \Bigr) \,\dd y\Biggr)^2,\\
F^{n,\eps,\delta,2}_r(v) &:=\frac{2}{\tau_n(r)^4}\sum_{k,m=1}^\infty \la_k^{r}\la_m^{r}\Biggl( \int_0^{\eps}  \Delta^n\ee_k(s)\Delta^n\ee_m(s)  \,\dd s \\
&\qquad\times \sum_{\ell=1}^{L(\delta)}\si^2(v-\eps,z^\delta_\ell) \int_{B^\delta_\ell}\phi_k(y)\phi_m(y)\,\dd y\Biggr)^2,\\
F^{n,\eps,\delta,3}_r(v) &:=\frac{2}{\tau_n(r)^4}\sum_{k,m=1}^\infty \la_k^{r}\la_m^{r}  \Biggl(\int_0^{\eps} \Delta^n\ee_k(s)\Delta^n\ee_m(s) \,\dd s \\
&\qquad\times\int_{D\setminus\bigcup_{\ell=1}^{L(\delta)} B^\delta_\ell} \phi_k(y)\phi_m(y)  \si^2(v-\eps,y) \,\dd y\Biggr)^2.
\end{align*}
Because $\si$ has uniformly bounded moments of all orders,
\[ \sup_{t\in[0,T]}\sup_{\ell =1,\dots,L(\delta)}\sup_{y\in B^\delta_\ell} \E[|\si^2(t,y) -\si^2(t,z^\delta_i) |^2]^{\frac12} \lec w_4(\delta;T) \]
with the notation from \eqref{L2-cont}. So following the calculations in the proof of Lemma~\ref{size-incr}, we deduce that for every $t\in[0,T]$,
\[  \lim_{\delta\to0}\limsup_{n\to\infty} \sup_{v\in[\eps,t]} \E[|F^{n,\eps,\delta,1}_r(v)|]= 0.\]
By a similar argument and the Cauchy--Schwarz inequality, 
\begin{align*} &\limsup_{n\to\infty} \sup_{v\in[\eps,t]}  \E[|F^{n,\eps,\delta,3}_r(v)|]\\
&\quad\lec \limsup_{n\to\infty} \Biggl( \frac{1}{\tau_n(r)^2} \sum_{k=1}^\infty \la_k^r \int_0^{\eps} (\Delta^n\ee_k(s))^2\,\dd s\int_{D\setminus\bigcup_{\ell=1}^{L(\delta)} B^\delta_\ell} \phi_k(y)^2\,\dd y \Biggr)^2\\
&\quad\lec \limsup_{n\to\infty}  \Biggl( \frac{1}{\tau_n(r)^2} \sum_{k=1}^\infty \la_k^{r-\ga}(1-\ee^{-\la_k^\ga\Del})  \int_{D\setminus\bigcup_{\ell=1}^{L(\delta)} B^\delta_\ell} \phi_k(y)^2\,\dd y \Biggr)^2. \end{align*}
By Lemma~\ref{convdens} (ii), this is $\lec |D\setminus\bigcup_{\ell=1}^{L(\delta)} B^\delta_\ell|^2$, which converges to $0$ as $\delta\to0$.

For $F^{n,\eps,\delta,2}_r(v)$, with small $\delta >0$ fixed, we consider the terms where $k=m$ and the terms where $k\neq m$ separately. By \eqref{mom} and \eqref{calc}, the sum of the terms where $k=m$ has a first moment that is bounded by
\beq\label{eq3}\begin{aligned}
	&\frac{2}{\tau_n(r)^4}\sum_{k=1}^\infty \la_k^{2r} \biggl( \int_0^\eps \Delta^n\ee_k(s)^2\,\dd s\biggr)^2\\
	&\quad\leq \frac{2}{\tau_n(r)^4}\sum_{k=1}^\infty \la_k^{2r-2\ga}  \Bigl( 1-\ee^{-\la_k^\ga\Del}  - \tfrac12\ee^{-2\la_k^\ga\eps} (\ee^{\la_k^\ga\Del}-1)^2 \Bigr)^2.
\end{aligned} \eeq
As in the proof of Lemma~\ref{size-incr}, the second part is of order $\tau_n(r)^{-4}\Del^4$, which goes to $0$ as $n \to \infty$. Regarding the first part, we have
\beq\label{sum}
\frac{1}{\tau_n(r)^4}\sum_{k=1}^\infty \la_k^{2r-2\ga} (1-\ee^{-\la_k^\ga\Del})^2\leq \Biggl(\sup_{k\in\N} \frac{\la_k^{r-\ga} (1-\ee^{-\la_k^\ga\Del})}{\tau_n(r)^2}\Biggr)\frac{1}{\tau_n(r)^2}\sum_{k=1}^\infty \la_k^{r-\ga} (1-\ee^{-\la_k^\ga\Del}).
\eeq
The term after the supremum converges to $K_r$ by Lemma~\ref{convdens} (i). At the same time, if $r\neq -\frac d2$, then $0< \frac1\ga(\ga-\frac d2-r)<1$. So choosing $\eps>0$ small enough, we obtain
\begin{align*}
\sup_{k\in\N} \frac{\la_k^{r-\ga} (1-\ee^{-\la_k^\ga\Del})}{\tau_n(r)^2} \leq \sup_{k\in\N} \frac{\la_k^{r-\ga}(\la_k^\ga\Del)^{(\ga-d/2-r+\eps)/\ga}}{\Del^{(\ga-d/2-r)/\ga}} = \Del^{\frac{\eps}{\ga}} \sup_{k\in\N}\la_k^{-\frac d2+\eps} = \la_1^{-\frac d2+\eps}\Del^{\frac \eps \ga} \to 0
\end{align*}
as $n\to\infty$. Similarly, if $r=-\frac d2$,
\begin{align*}
\sup_{k\in\N} \frac{\la_k^{-\frac d2-\ga} (1-\ee^{-\la_k^\ga\Del})}{\tau_n(-\frac d2)^2} \leq \sup_{k\in\N} \frac{\la_k^{-\frac d2-\ga} \la_k^\ga\Del}{\Del\lvert\log \Del\rvert} = \frac1{\lvert\log \Del\rvert} \sup_{k\in\N}\la_k^{-\frac d2} = \la_1^{-\frac d2}\frac1{\lvert\log \Del\rvert} \to 0.
\end{align*}
We conclude that the right-hand side of \eqref{sum} converges to $0$ as $n\to\infty$.

Finally, let us discuss the terms in $F^{n,\eps,\delta,2}_r(v)$ when $k\neq m$. Applying the Cauchy--Schwarz inequality to the $\dd s$-integral and realizing that, as before, the asymptotically relevant contribution of $\int_0^\eps (\Delta^n\ee_{j}(s))^2\,\dd s$ is $\la_{j}^{-\ga}(1-\ee^{-\la_{j}^\ga\Del})$, where $j \in \{k,m \}$, we only have to analyze further the following expression:
\begin{align*} &\frac{4}{\tau_n(r)^4}\sum_{k=1}^\infty\sum_{m=k+1}^\infty \la_k^{r-\ga}\la_m^{r-\ga}(1-\ee^{-\la_k^\ga\Del})(1-\ee^{-\la_m^\ga\Del})\\
&\qquad\times\Biggl(    \sum_{\ell=1}^{L(\delta)}\si^2(v-\eps,z^\delta_\ell) \int_{B^\delta_\ell}\phi_k(y)\phi_m(y)\,\dd y\Biggr)^2. \end{align*}
For any given $K\in\N$, we can further neglect all terms where $k\leq K$. Indeed, their sum is bounded in $L^1(\Om)$ by a constant times
\[ \frac{\Del}{\tau_n(r)^4}\sum_{k=1}^K \sum_{m=k+1}^\infty \la_k^r \la_m^{r-\ga}(1-\ee^{-\la_m^\ga\Del})\leq \frac{\Del}{\tau_n(r)^2}\sum_{k=1}^K \la_k^r \cdot \frac{1}{\tau_n(r)^2}\sum_{m=1}^\infty \la_m^{r-\ga}(1-\ee^{-\la_m^\ga\Del}).\]
We have $\Del/\tau_n(r)^2\to0$ and the second factor tends to $K_r$ by Lemma~\ref{convdens} (i). Thus, the previous display
converges to $0$ as $n\to\infty$, uniformly in $v\in[\eps,t]$.  As a consequence, only
\beq\label{ovF}\begin{aligned} \ov F^{n,\eps,\delta}_r(v)&:=\frac{4}{\tau_n(r)^4}\sum_{k=K+1}^\infty \sum_{m=k+1} \la_k^{r-\ga}\la_m^{r-\ga}(1-\ee^{-\la_k^\ga\Del})(1-\ee^{-\la_m^\ga\Del})\\
&\qquad\times\Biggl(    \sum_{\ell=1}^{L(\delta)}\si^2(v-\eps,z^\delta_\ell) \int_{B^\delta_\ell}\phi_k(y)\phi_m(y)\,\dd y\Biggr)^2 \end{aligned}\eeq
remains to be considered, where $K$ will be chosen shortly.

If $r>-\frac d2$, we decompose $\ov F^{n,\eps,\delta}_r(v)$ into $\ov F^{n,\eps,\delta,\theta,1}_r(v)+\ov F^{n,\eps,\delta,\theta,2}_r(v)$, where the first (resp., second) term comprises the indices $(k,m)$ satisfying $K<k<m\leq (1+\theta)k$ (resp., $k>K$ and $m>(1+\theta)k$) and where $\theta\in(0,1)$ is some fixed number for the moment. We can bound the $\dd y$-integral of $\phi_k\phi_m$ over $B^\delta_\ell$ in the previous display either by $1$ (using the Cauchy--Schwarz inequality) or by the estimate in Lemma~\ref{lem:innerprod} (recall that $\delta$ and the $B^\delta_\ell$ are fixed). Therefore, letting $\rho\in[0,1]$ be another parameter that we choose later, we obtain
\begin{align*} &\sup_{v\in[\eps,t]}\E[|\ov F^{n,\eps,\delta,\theta,2}_r(v)|] \\
&\quad\leq C_\delta \frac{4}{\tau_n(r)^4}\sum_{k=K+1}^\infty\sum_{m>(1+\theta)k}\la_k^{r-\ga}\la_m^{r-\ga}(1-\ee^{-\la_k^\ga\Del})(1-\ee^{-\la_m^\ga\Del})  \frac{\la_m^\rho \la_k^{\frac14 \rho}}{(\la_m-\la_k)^{2\rho}} \end{align*}
for some constant $C_\delta$ that is independent of $n$. Thanks to \eqref{eq:Weyl}, we can write $\la_k=C_D(1+\eta_k)k^{2/d}$ for all $k\in\N$ with a sequence  $(\eta_k)_{k\in\N}$   that converges to $0$. In particular, if $m>(1+\theta)k$, then
\[
  \la_m \geq \frac{1+\eta_m}{1+\eta_k} (1+\theta)^{\frac{2}{d}} \la_k,
  \qquad\text{therefore}\qquad  \la_m - \la_k \geq K(m,k,\theta) \, \la_k,
\]
where
\[
    K(m,k,\theta):= \frac{1+\eta_m}{1+\eta_k} (1+\theta)^{\frac{2}{d}} - 1 = \frac{1+\eta_m}{1+\eta_k} \left(\left[ (1+\theta)^{\frac{2}{d}} - 1\right] + 1 - \frac{1+\eta_k}{1+\eta_m}\right).
\]
Since $(1+\theta)^{\frac{2}{d}} - 1 \geq \frac{2\theta}{d}$, we have $K(m,k,\theta) \geq \frac{\theta}{2 d}$, for $m>k\geq K=K(\theta)\in \N$, where $K(\theta)$ is large enough so that
\[
   \frac{1+\eta_m}{1+\eta_k} \geq \frac12,\qquad
   1 - \frac{1+\eta_k}{1+\eta_m} \geq -\frac{\theta}{d}.
\]
Therefore,
 \begin{align*} &\sup_{v\in[\eps,t]}\E[|\ov F^{n,\eps,\delta,\theta,2}_r(v)|] \\
 &\quad\lec  \frac{1}{\tau_n(r)^4}\sum_{k=K+1}^\infty\sum_{m>(1+\theta)k}\la_k^{r-\ga}\la_m^{r-\ga}(1-\ee^{-\la_k^\ga\Del})(1-\ee^{-\la_m^\ga\Del})   \la_m^\rho \la_k^{\frac14 \rho}\la_k^{-2 \rho}\\
&\quad\lec \frac{1}{\tau_n(r)^4}\sum_{k=1}^\infty \la_k^{r-\frac74 \rho-\ga} (1-\ee^{-\la_k^\ga\Del}) \sum_{m=1}^\infty \la_m^{ r-\ga+\rho}(1-\ee^{-\la_m^\ga\Del}). \end{align*}
Let $\rho>0$ be small enough such that $r_1=r-\frac 74 \rho>-\frac d2$ and $r_2=r+\rho<\ga-\frac d2$. Then, by Lemma~\ref{convdens} (i), the sum over $k$ is of order $\tau_n(r_1)^2$, while the sum over $m$ is of order $\tau_n(r_2)^2$. As a consequence,
the formula \eqref{taunr} for $\tau_n(r)$ shows
\beq\label{F2-2}   \lim_{n\to\infty} \sup_{v\in[\eps,t]}\E[|\ov F^{n,\eps,\delta,\theta,2}_r(v)|] \lec \lim_{n\to\infty} \frac{\tau_n(r_1)^2\tau_n(r_2)^2}{\tau_n(r)^4} = \lim_{n\to\infty} \Del^{\frac1\ga(\frac74\rho-\rho)}= 0. \eeq

Next, we consider $\ov F^{n,\eps,\delta,\theta,1}_r(v)$. Bounding the second moment of the sum in parentheses in \eqref{ovF} by a constant and using $1-\ee^{-x}\leq x^{\frac1\ga(\ga-\frac d2-r)}$, we obtain
\beq\label{F1}\begin{aligned}
\sup_{v\in[\eps,t]} \E[|\ov F^{n,\eps,\delta,\theta,1}_r(v)|] &\lec \frac{1}{\tau_n(r)^4}\sum_{k=1}^\infty \sum_{k<m\leq (1+\theta)k}  \la_k^{r-\ga}\la_m^{r-\ga}(1-\ee^{-\la_k^\ga\Del})(1-\ee^{-\la_m^\ga\Del})\\
&\lec \frac{1}{\tau_n(r)^2}\sum_{k=1}^\infty \la_k^{r-\ga} (1-\ee^{-\la_k^\ga\Del})  \sum_{k<m\leq (1+\theta)k} m^{-1}\\
&\leq \frac{1}{\tau_n(r)^2}\sum_{k=1}^\infty \la_k^{r-\ga} (1-\ee^{-\la_k^\ga\Del}) (\log((1+\theta)k)-\log k)\\
&=\log(1+\theta)\frac{1}{\tau_n(r)^2}\sum_{k=1}^\infty \la_k^{r-\ga} (1-\ee^{-\la_k^\ga\Del})  \\
&\to K_r \log(1+\theta)
\end{aligned}\eeq
by Lemma~\ref{convdens} (i). As a consequence,
\[ \lim_{\theta\to0}\limsup_{n\to\infty}\sup_{v\in[\eps,t]} \E[|\ov F^{n,\eps,\delta,\theta,1}_r(v)|]  =0,  \]
which completes the proof of \eqref{conv-var} for $r>-\frac d2$.

If $r=-\frac d2$, we go back to \eqref{ovF} and let $\ov F^{n,\eps,\delta,\theta,1}_{-d/2}(v)$ (resp., $\ov F^{n,\eps,\delta,\theta,2}_{-d/2}(v)$) contain the indices $(k,m)$ with $K<k<m\leq k^\theta$ (resp., $k>K$ and $m>k^\theta$), where $\theta>1$ and $K=K(\theta)>0$ is chosen such that $|\eta_k|\leq \frac12$ for all $k>K$ and $m^{\frac2d}\geq 6 m^{\frac 2{d\theta}}$ for all $m>K^\theta$. Then $k>K$ and $m>k^\theta$ imply
\begin{align*} \la_m-\la_k &= C_D(1+\eta_m)m^{\frac 2d} - C_D(1+\eta_k)k^{\frac 2d} \geq \tfrac12C_Dm^{\frac 2d} -\tfrac32 C_Dm^{\frac2{d\theta}}\geq\tfrac14 C_Dm^{\frac 2d}\\
&\begin{cases} \gtrsim \la_m,  \\  \geq\frac14 C_Dk^{\frac 2d \theta}\gtrsim \la_k^\theta. \end{cases}\end{align*}

So arguing as above (with $\rho=1$) and bounding $(\la_m-\la_k)^2\gtrsim \la_m\la_k^\theta$, we derive
 \begin{align*} &\sup_{v\in[\eps,t]}\E[|\ov F^{n,\eps,\delta,\theta,2}_{-d/2}(v)|]\\
  &\quad\lec  \frac{1}{\tau_n(-\frac d2)^4}\sum_{k=K+1}^\infty\sum_{m>k^\theta}\la_k^{-\frac d2-\ga}\la_m^{-\frac d2-\ga}(1-\ee^{-\la_k^\ga\Del})(1-\ee^{-\la_m^\ga\Del})     \la_k^{\frac14 -\theta} \\
&\quad\lec  \frac{1}{\tau_n(-\frac d2)^4}\sum_{k= 1}^\infty\la_k^{-\frac d2+\frac14 -\theta-\ga}(1-\ee^{-\la_k^\ga\Del})\sum_{m=1}^\infty\la_m^{-\frac d2-\ga}(1-\ee^{-\la_m^\ga\Del}).
\end{align*}
By Lemma \ref{convdens}, the sum over $m$ is $\calo(\tau_n(-\frac d2)^2)$, whereas, as seen in the proof of Lemma \ref{size-incr}, the sum over $k$ is $\calo(\tau_n(-\frac d2+\frac14-\theta))=\calo(\Del)$. Therefore,
 \begin{align*} \lim_{n\to\infty}\sup_{v\in[\eps,t]}\E[|\ov F^{n,\eps,\delta,\theta,2}_{-d/2}(v)|] &\lec\lim_{n\to\infty} \frac{\Del}{\tau_n(-\frac d2)^2} = \lim_{n\to\infty} \frac1{\lvert\log \Del\rvert} = 0. \end{align*}

Analogously to \eqref{F1}, we further have
\begin{align*} &\sup_{v\in[\eps,t]} \E[|\ov F^{n,\eps,\delta,\theta,1}_{-d/2}(v)|] \\
&\quad\lec \frac{\theta-1}{\tau_n(-\frac d2)^2\lvert\log\Del\rvert}\sum_{k=1}^\infty \la_k^{-\frac d2-\ga} (1-\ee^{-\la_k^\ga\Del})\log k\\
&\quad\lec \frac{\theta-1}{\lvert\log\Del\rvert^2}\sum_{k=1}^{[\Del^{-d/(2\ga)}]} k^{-1} \log k +\frac{\theta-1}{\Del\lvert\log\Del\rvert^2}\sum_{k=[\Del^{-d/(2\ga)}]+1}^\infty k^{-1-\frac2d\ga}  \log k. \end{align*}
Bounding $k^{-\frac 2d \ga}\leq(1+[\Del^{-d/(2\ga)}])^{-\frac 2d \ga}\leq \Del$, we conclude that
\begin{align*} \limsup_{n\to\infty} \sup_{v\in[\eps,t]} \E[|\ov F^{n,\eps,\delta,\theta,1}_{-d/2}(v)|] &\leq \lim_{n\to\infty} \frac{2(\theta-1)}{\lvert\log\Del\rvert^2}\log \Del^{-\frac{d}{2\ga}}(1+\log \Del^{-\frac{d}{2\ga}})\lec \theta-1. \end{align*}
The proof of \eqref{conv-var} when $r=-\frac d2$ is completed by sending $\theta\downarrow 1$.
\epr

\blem\label{LLN-condexp}
For every $t>0$,
\[
\lim_{\eps\to0} \limsup_{n\to\infty} \bbe\Biggl[\Biggl\lvert \Del\sumt  f\biggl( \frac{\wt h^{n,i,\eps}_r}{\tau_n(r)}\biggr) -V^r_f(u,t)\Biggr\rvert\Biggr] = 0.\]
\elem
\bpr
We write the difference in the previous line as a sum of four terms: 
\begin{align*}
G^{n,\eps,1}_r &:= \Del \sum_{i=[\eps/\Del]+1}^{[t/\Del]} \Biggl\{  f\Biggl(\frac{\wt h^{n,i,\eps}_r}{\tau_n(r)}\Biggr)  - f\Biggl( \sqrt{\frac{K_r}{|D|} \int_D\si^2(i\Del-\eps,y)\,\dd y}\Biggr) \Biggr\},\\
G^{n,\eps,2}_r &:=  \int_{\Del[\eps/\Del]}^{\Del[t/\Del]}  \sum_{i=[\eps/\Del]+1}^{[t/\Del]} \Biggl\{ f\Biggl( \sqrt{\frac{K_r}{|D|} \int_D \si^2(i\Del-\eps,y)\,\dd y}\Biggr)\\
&\qquad- f\Biggl( \sqrt{\frac{K_r}{|D|} \int_D \si^2(v-\eps,y)\,\dd y}\Biggr) \Biggr\}\bone_{[(i-1)\Del,i\Del)}(v) \,\dd v,\\
G^{n,\eps,3}_r &:= \int_{\Del[\eps/\Del]}^{\Del[t/\Del]} \Biggl\{  f\Biggl( \sqrt{\frac{K_r}{|D|} \int_D\si^2(v-\eps,y)\,\dd y}\Biggr)- f\Biggl( \sqrt{\frac{K_r}{|D|} \int_D \si^2(v,y)\,\dd y}\Biggr)\Biggr\}\,\dd v,\\
G^{n,\eps,4}_r &:=\Del\sum_{i=1}^{[\eps/\Del]}  f\Biggl( \frac{\wt h^{n,i,\eps}_r}{\tau_n(r)}\Biggr) -\int_0^{\Del[\eps/\Del]} f\Biggl( \sqrt{\frac{K_r}{|D|} \int_D\si^2(v,y)\,\dd y}\Biggr)\,\dd v\\
&\qquad-\int_{\Del[t/\Del]}^{t} f\Biggl( \sqrt{\frac{K_r}{|D|} \int_D \si^2(v,y)\,\dd y}\Biggr)\,\dd v.
\end{align*}

Employing the estimate \eqref{fdiff} and Lemma~\ref{convP} and arguing as in the beginning of the proof of Lemma~\ref{LLN-discr} to get \eqref{toshow1} for $\wt h^{n,i,\eps}_r$, we can show that $\bbe[|G^{n,\eps,1}_r|]\to0$ as $n\to\infty$. Using the $L^2(\Om)$-continuity of $\si$ in addition (see Remark~\ref{L2-rem}), it follows similarly to the proof of Lemma~\ref{LLN-discr}, that $\bbe[|G^{n,\eps,2}_r|]\to0$ as $n\to\infty$ and $\limsup_{n\to\infty} \bbe[|G^{n,\eps,3}_r|] \to 0$ as $\eps\to0$. Furthermore, as the reader may easily verify by using \eqref{toshow2} for $\wt h^{n,i,\eps}_r$ and \eqref{mom},  also $\lim_{\eps\to0}\limsup_{n\to\infty} \bbe[|G^{n,\eps,4}_r|]= 0$.
\epr

\bpr[Proof of Corollary~\ref{Holder}] With the same arguments as for \eqref{toshow-2} and \eqref{toshow2}, one can show that for all $p,T>0$, there is a constant $C=C(p,r,\ga,\si,T)>0$ such that for all $s,t\in[0,T]$
\[ \E[\|u(t,\cdot)-u(s,\cdot)\|^p_{H_r}]^{\frac1p} \leq C|t-s|^{\al(r)}  \]
($\leq C|t-s|^{\al(r)-\eps}$ when $r=-\frac d2$, where $\eps$ can be taken arbitrarily small). So the first part of the corollary follows from Kolmogorov's continuity theorem (see \cite[Theorem~3.3]{DaPrato14}). For the second part, if $\si\equiv1$, observe from Lemma~\ref{size-incr} that for all $0<T_1<T_2<\infty$, there is $C'=C'(r,\ga,T_1,T_2)$ such that for all $s,t\in[T_1,T_2]$,
\[ \E[\|u(t,\cdot)-u(s,\cdot)\|^2_{H_r}]\geq C'|t-s|^{2\al(r)}.  \]
Hence, the second assertion  can be obtained with the same arguments as in \cite[Theorem~4]{Dalang05}.
\epr


\subsection{Comparison between Theorem~\ref{LLN}  and Theorem~\ref{LLN2}}\label{compare}

We conclude this paper by comparing Theorem~\ref{LLN} (where $r<-\frac d2$) and Theorem~\ref{LLN2} (where $-\frac d2 \leq r<\ga-\frac d2$) in two remarks.
\brem\label{rem-tight} Let us explain why in Theorem~\ref{LLN2}, we did not study functionals as general as in Theorem~\ref{LLN}. As Lemma~\ref{tight} shows, the family
\[ \calu' := \Biggl\{ \frac{u(i\Del,\cdot)-u((i-1)\Del,\cdot)}{\tau_n(r)}: \quad n\in\N,~i=1,\dots,[T/\Del]\Biggr\} \]
is tight for $r<-\frac d2$ and $\tau_n(r)=\sqrt{\Del}$, which makes it possible to prove Lemma \ref{remove-drift}.

By contrast, if $-\frac d2\leq r<\ga-\frac d2$ and $\si$ is, say, identically $1$, then there is no normalizing sequence $\tau_n(r)$ such that $\calu'$ becomes tight in $H_r$ without all limit points being zero. Indeed, 
it follows from Lemma~\ref{size-incr} that for $-\tfrac d2 \leq r < \ga-\tfrac d2$, the normalization $\tau_n(r)$ must be chosen as in \eqref{taunr} (or at least, of the same order) to avoid degenerate limits. Since $\calu'$ is a collection of centered Gaussian vectors on $H_r$, we know from \cite[Example~3.8.13~(iv)]{Bogachev98} that if $\calu'$ were tight, then the series $\sum_{k=1}^\infty \bbe[\langle \Delta^n_i u/\tau_n(r), b_k\rangle_{H_r}^2]$ would converge uniformly in $n$ and $i$, where $b_k$ are the orthonormal basis functions from \eqref{basis_b} and $\Delta^n_i u= u(i\Del,\cdot)-u((i-1)\Del,\cdot)$. In other words, we would have
\beq\label{false} \lim_{m\to\infty} \sup_{n\in\N}\sup_{i=1,\ldots,[T/\Del]}\sum_{k=m}^\infty \bbe\Biggl[\biggl\langle \frac{\Delta^n_i u}{\tau_n(r)}, b_k\biggr\rangle_{H_r}^2\Biggr] = 0. \eeq
But this is not true because from the calculations in \eqref{series-calc}, we deduce that
\begin{align*}&\sum_{k=m}^\infty \bbe\Biggl[\biggl\langle \frac{\Delta^n_i u}{\tau_n(r)}, b_k\biggr\rangle_{H_r}^2\Biggr]\\
&\quad= \frac{1}{\tau_n(r)^2}\Biggl(\sum_{k=m}^\infty \la_k^{r-\ga}(1-\ee^{-\la_k^\ga\Del})  - \frac12\sum_{k=m}^\infty \la_k^{r-\ga} \ee^{-2\la_k^\ga i\Del}  (\ee^{\la_k^\ga\Del}-1)^2\Biggr). \end{align*}
The expression within the parentheses is exactly the last line of \eqref{calc}, except that summation starts at $k=m$. By \eqref{second-term}, if $i\Del\geq\eps$ for some $\eps>0$, then the second part, divided by $\tau_n(r)^2$, vanishes as $n\to\infty$. In the notation of the proof of Lemma~\ref{size-incr}, the first part equals
\beq\label{eq4} \frac{1}{\tau_n(r)^2}A^{n,1}_r- \frac{1}{\tau_n(r)^2}\sum_{k=1}^{m-1}\la_k^{r-\ga}(1-\ee^{-\la_k^\ga\Del}). \eeq
For any fixed value of $m\in\N$, the second term is bounded by $(\sum_{k=1}^{m-1} \la_k^r)\Del/\tau_n(r)^2$, which converges to $0$ as $n\to\infty$. At the same time, the first term converges to the constant $K_r$ of \eqref{Kr}, as seen in the proof of the lemma. This is a nonzero limit that is independent of the value of $m$, so \eqref{false} must be false and therefore $\calu'$ cannot be tight.
\erem

\brem\label{difference} There are further important differences between Theorem~\ref{LLN} and Theorem~\ref{LLN2}. Let us take $\si\equiv1$  and consider quadratic variations (i.e., $p=2$ in \eqref{Vnp}, so $f(x) = x^2$) in the following to simplify the discussion.
\benu
	\item If $r<-\frac d2$, Theorem~\ref{LLN} gives a genuine law of large numbers in the following sense: As $n$ tends to infinity, most of the squared normalized increments $\|\Delta^n_i u\|^2_{H_r}/\tau_n(r)^2$ (where ``most'' means that  $i\Del\geq\eps$ for some fixed but arbitrary $\eps>0$) have equal mean (namely $K_r$ by Lemma~\ref{size-incr}) and bounded variance (see \eqref{toshow-2}). Moreover, they are only weakly correlated to each other, so a classical $L^2(\Om)$-argument gives the law of large numbers; cf.\ Lemma~\ref{LLN-core-1}.
	\item By contrast, if $-\frac d2\leq r<\ga-\frac d2$, Theorem~\ref{LLN2} is  a \emph{degenerate} law of large numbers. Here, not only do we have $\bbe[\|\Delta^n_i u\|_{H_r}^2]/\tau_n(r)^2 \to K_r$ (i.e., convergence of means) for most values of $i$, but we actually have $\bbe[|\|\Delta^n_i u\|_{H_r}^2/\tau_n(r)^2 - K_r|]$ (i.e., convergence in mean); cf.\ Lemma~\ref{convP}. In the case $r<-\frac d2$, only the \emph{sum} $\Del\sum_{i=1}^{[t/\Del]} \|\Delta^n_i u\|_{H_r}^2/\tau_n(r)^2$ converges in mean, not the individual summands. In particular, no averaging or cancellation argument as in Lemma~\ref{LLN-core-1}, which would be typical of a genuine law of large numbers, is needed for $-\frac d2\leq r<\ga-\frac d2$. This also explains why if we compare \eqref{2pvar-add} and \eqref{p-var-add}, then only in \eqref{2pvar-add} is there a factor (namely $2^pB_p(x_1,\ldots,x_p)$) that is related to moment formulas.
\eenu
\erem

\begin{appendix}
	\setcounter{equation}{0}
\setcounter{theorem}{0}
\section*{}\label{appn}
We start with some analytic preliminaries needed for the proof of Proposition~\ref{exist}.
\blem\label{A1} Let $D$ be an open bounded subset of $\R^d$. Then the eigenvalues and eigenfunctions of $-\Delta$ have the following properties:
\benu
\item The eigenvalues $(\la_n)_{n\in\N}$ satisfy \emph{Weyl's law}, that is,
\beq\label{eq:Weyl} \lim_{n\to\infty} \frac{\la_n}{n^{2/d}} =C_D:=\frac{4\pi \Ga(1+\frac d2)^{2/d}}{|D|^{2/d}}. \eeq
\item For all $n\in\N$,
\beq\label{eq:Linfty} \|\phi_n\|_{L^\infty(D)} \leq \biggl(\frac{2\ee \la_n}{\pi d}\biggr)^{d/4}. \eeq
\item For every $x\in D$,  the following asymptotics are valid as $\la\to\infty$ (where $f(x)\sim g(x)$ means $f(x)/g(x)\to1$):
\bit
\item for $r<-\frac d 2$,
\beq\label{est1} \sum_{k:\,\la_k>\la} \la_k^r |\phi_k(x)|^{2} \sim -\frac{1}{(4\pi)^{d/2}\Ga(1+\frac d2)(1+\frac 2d r)} \la^{\frac d2+r}; \eeq
\item for $r>-\frac d2$,
\beq\label{est2} \sum_{k:\,\la_k\leq\la} \la_k^r |\phi_k(x)|^{2} \sim \frac{1}{(4\pi)^{d/2}\Ga(1+\frac d2)(1+\frac 2d r)} \la^{\frac d2+r}; \eeq
\item for $r=-\frac d2$,
\beq\label{est3} \sum_{k:\,\la_k\leq\la} \la_k^r |\phi_k(x)|^{2} \sim \frac{1}{(4\pi)^{d/2}\Ga(\frac d2)}\log \la. \eeq
\eit
\item If $D$ additionally has the cone property, then, for any $r<-\frac d2$, there is $\eps>0$ such that
\beq\label{summable} \sup_{x\in D} \sum_{k=1}^\infty \la_k^{r} |\phi_k(x)|^{2+\eps}<\infty. \eeq
\eenu
\elem
\bpr Part (i) is well known; see \cite[Theorem~1.11]{Birman80} or \cite{Geisinger14}. Part (ii) is shown in \cite[Lemma~3.1]{Davies74}. For (iii), we use \cite[Equation~(0.7)]{Garding53}, which shows, for every $x\in D$, that
\[ V(t):=V(t;x):=\sum_{1\leq k\leq t} |\phi_k(x)|^2 \sim \frac{t}{|D|} \]
as $t\to\infty$.
Therefore, by Abel's summation formula (see \cite[Theorem~4.2]{Apostol76}),
\begin{align*} \sum_{1\leq k\leq t} k^{\frac 2d r} |\phi_k(x)|^2 &= t^{\frac 2d r}V(t) -\frac{2}{d}r \int_1^t s^{\frac2d r-1} V(s)\,\dd s \sim  \frac{t^{\frac2d r+1}}{|D|} - \frac{2r}{d\lvert D\rvert}\int_1^t s^{\frac2d r}\,\dd s\\
&\sim \begin{cases} \frac{1}{(\frac 2d r+1)\lvert D\rvert}t^{\frac 2d r+1}&\text{if } r>-\frac d2, \\ \frac{1}{|D|}\log t &\text{if } r=-\frac d2.\end{cases} \end{align*}
As $\la_k \sim C_Dk^{2/d}$ by \eqref{eq:Weyl}, we derive
\[ \sum_{k:\, \la_k\leq\la} \la_k^r |\phi_k(x)|^{2} \sim C_D^r \sum_{1\leq k\leq C_D^{-d/2}\la^{d/2}} k^{\frac 2d r}|\phi_k(x)|^2 \sim \begin{cases} \frac{C_D^{- d/2}}{(\frac 2d r+1)\lvert D\rvert}\la^{r+\frac d2}&\text{if } r>-\frac d2, \\ \frac{dC_D^{-d/2}}{2|D|}\log \la &\text{if } r=-\frac d2,\end{cases} \]
as $\la\to\infty$, which becomes \eqref{est2} and \eqref{est3} after simplifying the constants. If $r<-\frac d2$, we obtain, with similar arguments,
\begin{align*} \sum_{k:\, k> t} k^{\frac 2d r} |\phi_k(x)|^2 &= \lim_{K\to\infty} \Biggl( K^{\frac 2d r}V(K) - t^{\frac2d r} V(t) -\frac{2}{d}r \int_t^K s^{\frac2d r-1} V(s)\,\dd s\Biggr)\\
&= - t^{\frac2d r} V(t) -\frac{2}{d}r \int_t^\infty s^{\frac2d r-1} V(s)\,\dd s\\
&\sim  -\frac{t^{\frac2d r+1}}{|D|} - \frac{2r}{d\lvert D\rvert}\int_t^\infty s^{\frac2d r}\,\dd s= - \frac{1}{(\frac 2d r+1)\lvert D\rvert}t^{\frac 2d r+1} \end{align*}
and consequently,
\[ \sum_{\la_k>\la} \la_k^r |\phi_k(x)|^{2} \sim C_D^r \sum_{k> C_D^{-d/2}\la^{d/2}} k^{\frac 2d r}|\phi_k(x)|^2 \sim  - \frac{C_D^{- d/2}}{(\frac 2d r+1)\lvert D\rvert}\la^{r+\frac d2},   \]
which is \eqref{est1}.

For (iv), we use (ii), \cite[Theorem~8.2]{Agmon65}, and the hypothesis $r <-\frac d2$ to deduce
\begin{align*}
\sum_{\la_k\geq\la} \la_k^{r} |\phi_k(x)|^{2+\eps} &=\sum_{k=0}^\infty \sum_{\la_k\in [2^n\la, 2^{n+1}\la)} \la_k^{r} |\phi_k(x)|^{2+\eps} \lec \sum_{n=0}^\infty \sum_{\la_k\in [2^n\la, 2^{n+1}\la)} \la_k^{r+ \frac d 4 \eps} \phi_k(x)^2 \\
&\lec \la^{r+ \frac d 4 \eps}\sum_{n=0}^\infty 2^{nr} \sum_{\la_k\leq 2^{n+1}\la}\phi_k(x)^2 \lec \la^{r+ \frac d 4 \eps} \sum_{n=0}^\infty 2^{nr} (2^{n+1}\la)^{\frac d2} \lec \la^{\frac d2 +r+ \frac d 4 \eps},
\end{align*}
which converges to $0$ as $\la\to\infty$, uniformly in $x\in D$, if $\eps$ is small enough. This proves \eqref{summable}.
\epr

\bpr[Proof of Proposition~\ref{exist}] The existence of a jointly measurable and adapted solution of \eqref{SHE} will be a consequence of Theorem 4.2.1 in \cite{Dalang20} once we show that
\beq\label{L2kernel}   \int_0^t \dd s \, \sup_{x\in\R^d} \int_D \dd y\, g^2(s;x,y) < \infty.
\eeq

According to \cite[Theorem 8.2]{Agmon65}, because $D$ satisfies the cone property, there is a finite constant $C>0$ such that for all $x \in D$ and all $t \geq 0$,
\beq\label{eA.8}
   \wt V(t,x) := \sum_{k:\, \lambda_k \leq t} |\phi_k(x)|^2 \leq C t^{d/2}.
\eeq

Using the orthonormal property of $(\phi_k)_{k\in\N}$, we see that
\begin{align*} \int_D g(s;x,y)^2\,\dd y 
&= \sum_{k=1}^\infty \phi_k(x)^2\ee^{-2\la_k^\ga s}
= \int_{\lambda_1-}^{+\infty} \ee^{-2u^\gamma s} \, \wt V(\dd u,x).
\end{align*}
By Abel's summation formula (see \cite[Theorem~4.2]{Apostol76}), this is equal to
\[
   - \int_{\lambda_1}^{+\infty} \wt V(u,x) (-2\gamma u^{\gamma-1} s) \ee^{-2u^\gamma s}\, \dd u \leq 2\gamma C \int_{\lambda_1}^{+\infty} u^{\frac{d}{2} + \gamma -1} s \ee^{2 u^\gamma s} \, \dd u.
\]
Therefore,
\[
   \int_0^t \dd s \, \sup_{x\in\R^d} \int_D \dd y\, g^2(s;x,y) \leq 2\gamma C \int_{\lambda_1}^{+\infty} \dd u\, u^{\frac{d}{2} + \gamma -1}\int_0^t \dd s \, s\, \ee^{2 u^\gamma s}.
\]
Since $\int_0^\infty se^{-as}\dd s = a^{-2}$ ($a>0$), this is bounded above by
\[
   2\gamma C \int_{\lambda_1-}^{+\infty} u^{\frac{d}{2} + \gamma -1} (2 u^\gamma)^{-2}\, \dd u = \tfrac{\gamma C}{2} \int_{\lambda_1-}^{+\infty}  u^{\frac{d}{2} - \gamma -1}\, \dd u < \infty
\]
because $\tfrac{d}{2} - \gamma <0$. This implies \eqref{L2kernel}.


In order to establish the $L^p$-continuity of $u$, according to \cite[Section 4.2.2]{Dalang20}, it suffices to show that
\begin{align}\label{L2cont} 
\lim_{(\tau,h)\to0}  \int_0^t \int_D (g(t+\tau-s;x+h,y)-g(t-s;x,y))^2\,\dd y\,\dd s &=0,\\
\lim_{\tau \downarrow 0}  \int_t^{t+\tau}\int_D g^2(t+\tau-s;x,y)\,\dd y\,\dd s &= 0, \label{L2cont2}\end{align}
for all $t>0$ and $x\in D$. Once $L^2$-continuity is established, predictability follows immediately \cite[Proposition 2]{Dalang98}.

In order to establish \eqref{L2cont}, because $\int_0^t (\ee^{-\la_k^\ga (s+\tau)}-\ee^{-\la_k^\ga s})^2\,\dd s = (1-\ee^{-\la_k^\ga \tau})^2(1-\ee^{-2\la_k^\ga t})/(2\la_k^\ga)$, one easily checks that the integral in \eqref{L2cont} is bounded by
\begin{align*} &2\int_0^t\int_D (g(s+\tau;x+h,y)-g(s;x+h,y))^2\,\dd y\,\dd s \\
&\quad\quad+2\int_0^t\int_D (g(s;x+h,y)-g(s;x,y))^2\,\dd y\,\dd s\\
&\quad \leq \sum_{k=1}^\infty \frac{(1-\ee^{-\la_k^\ga \tau})^2}{\la_k^\ga}\phi_k(x+h)^2 + \sum_{k=1}^\infty \frac{(\phi_k(x+h)-\phi_k(x))^2}{\la_k^\ga}. \end{align*}
As $(\tau,h)\to0$, the terms within both sums converge to $0$ for each $k\in\N$. Moreover, $\sum_{k=1}^\infty \la_k^{-\ga}\delta_k$ is a finite measure on $\N$ by Lemma~\ref{A1} (i) and the hypothesis $\ga>\frac d2$. So by Proposition~4.12 and the discussion before Lemma~4.10 in \cite{Kallenberg02}, \eqref{L2cont} is implied by $\sup_{x\in D} \sum_{k=1}^\infty (2\la_k^\ga)^{-1}|\phi_k(x)|^{2+\eps}<\infty$, which was shown in \eqref{summable} for some small $\eps>0$.

Concerning \eqref{L2cont2}, the same arguments show that
\[ \int_t^{t+\tau}\int_D g^2(t+\tau-s;x,y)\,\dd y\,\dd s = \int_0^\tau \int_D g^2(s;x,y)\,\dd y\,\dd s = \sum_{k=1}^\infty \frac{1-\ee^{-\la_k^\ga \tau}}{2\la_k^\ga} \phi_k(x)^2 \to 0 \]
as $\tau \downarrow 0$.
\epr

\bpr[Proof of Proposition~\ref{uinHr}] Observe from \eqref{akt} and the orthogonality of $(\phi_k)_{k\in\N}$ that $(a_k(t))_{k\in\N}$ is a sequence of independent centered Gaussian random variables. Thus, by \cite[Theorem~4.17]{Kallenberg02} and \cite[Theorem~22.7]{Billingsley95}, the almost-sure finiteness of $\|u(t,\cdot)\|^2_{H_r} = \sum_{k=1}^\infty \la_k^r (a_k(t))^2$ is equivalent to the convergence in probability of $\sum_{k=1}^n \la_k^{r/2} a_k(t)$ as $n\to\infty$. By Gaussianity, this is in turn equivalent to the summability of the variances, that is, to $r<\ga-\frac d2$ by \eqref{eq:L2finite}.
\epr

\blem\label{convdens} Suppose that $D\subseteq\R^d$ is open and bounded and that $-\frac d2\leq r<\ga-\frac d2$. Let
\begin{align*}
H_n&:=\frac1{\tau_n(r)^2}   \sum_{k=1}^\infty \la_k^{r-\ga} (1-\ee^{-\la_k^\ga\Del}),\\
H_n(y)&:= \frac1{\tau_n(r)^2}   \sum_{k=1}^\infty \la_k^{r-\ga} (1-\ee^{-\la_k^\ga\Del})\phi_k(y)^2,\qquad y\in D, \end{align*}
with $\tau_n(r)$ from \eqref{taunr}. Further recall the definition of $K_r$ from \eqref{Kr}.
\benu
\item As $n\to\infty$, we have $H_n\to K_r$ and $H_n(y)\to K_r/|D|$ for every $y\in D$.
\item If $D$ additionally satisfies the cone property and $D_0$ is a measurable subset of $D$, then
\beq\label{limsup} \limsup_{n\to\infty} \int_{D_0} H_n(y)\,\dd y \leq C|D_0|\eeq
with a constant $C\in(0,\infty)$ that does not depend on $D_0$.
\eenu
\elem
\bpr Consider (i) first. We only show the claim for $H_n(y)$; the assertion for $H_n$ follows analogously once we replace $\phi_k(y)^2$ by $1$ in the subsequent arguments (in this case, there would be an additional factor $|D|$ in \eqref{eA.12}, which explains that the limit of $H_n$ differs from the limit of $H_n(y)$ also by a factor of $|D|$). If $r\neq -\frac d2$, we fix $y\in D$ and
let
\begin{align*} U(\la)&:=U(\la;y):=\sum_{k:\,\la_k\leq \la} \phi_k(y)^2,\qquad \la>0,\\
k(z)&:=z^{\ga-r}(1-\ee^{-z^{-\ga}}),\qquad z>0,\\
M(z):=M(z;y)&:= \int_0^\infty k(z/\la)\, \dd U(\la) = z^{\ga-r}\sum_{k=1}^\infty \la_k^{r-\ga}(1-\ee^{-\la_k^\ga/z^{\ga}})\phi_k(y)^2 ,\qquad z>0.
\end{align*}
so that $M(z)$ is the Mellin--Stieltjes transform of $U$ \cite[(4.0.3)]{Bingham87}. Then
\beq\label{Hn} H_n(y)=\Del^{-\frac1\ga(\ga-\frac d2 -r)}\sum_{k=1}^\infty \la_k^{r-\ga} (1-\ee^{-\la_k^\ga\Del})\phi_k(y)^2 = \Del^{\frac {d}{2\ga}}M(\Del^{-\frac1\ga}). \eeq

We use an Abelian theorem from \cite{Bingham87} to determine the behavior of $M(z)$ as $z\to\infty$. To this end, observe from \cite[Equation~(0.7)]{Garding53} and \eqref{eq:Weyl} that
\beq\label{eA.12} U(\la)\sim \sum_{1\leq k\leq C_D^{-d/2}\la^{d/2}} \phi_k(y)^2\sim (2\pi)^{-d}\frac{\pi^{\frac d2}}{\Ga(1+\frac d2)}\la^{\frac d2},\qquad \la\to\infty. \eeq
 Moreover, since $\la_1>0$, $U(\la)$ is identically $0$ for $0<\la<\la_1$. Next, pick $\al,\beta>0$ such that $(-r)\vee 0<\al<\frac d2 <\beta<\ga -r$, which is possible because $-\frac d 2<r<\ga-\frac d2$ by assumption. Because $z\mapsto k(z)$ is continuous on $(0,\infty)$, increasing at $0$ and eventually decreasing for large $z$, there is $n_0\in\N$ such that
\begin{align*} &\sum_{n\in\Z} (\ee^{-\al n}\vee\ee^{-\beta n})\sup_{\ee^n \leq z\leq \ee^{n+1}} k(z)\\
&\quad \leq \sum_{n= -\infty}^{-n_0-1} \ee^{-\beta n}\ee^{(n+1)(\ga-r)} + \sum_{n=-n_0}^{n_0} (\ee^{-\al n}\vee\ee^{-\beta n}) \sup_{z\in(0,\infty)}k(z)+\sum_{n=n_0+1}^\infty \ee^{-\al n}\ee^{-nr}\\
&\quad<\infty.\end{align*}
We conclude from \cite[Theorem~4.4.2]{Bingham87} that
\[ M(z)\sim \frac d2(2\pi)^{-d}\frac{\pi^{\frac d2}}{\Ga(1+\frac d2)}\check k(\tfrac d2)z^{\frac d2},\qquad z\to\infty, \]
where, by a substitution $z^\ga\mapsto y$ and integration by parts,
\begin{align*} \check k(\tfrac d2)&:=\int_0^\infty z^{\frac d2-1}k(z^{-1})\,\dd z = \int_0^\infty z^{\frac d2-1}z^{r-\ga}(1-\ee^{-z^\ga})\,\dd z=\frac1\ga \int_0^\infty y^{\frac{r}{\ga}-2+\frac {d}{2\ga}}(1-\ee^{-y}) \,\dd y \\
&=\frac{-1}{r-\ga+\frac {d}{2}}\int_0^\infty  y^{\frac{r}{\ga}-1+\frac {d}{2\ga}}\ee^{-y}\,\dd y=\frac{-1}{r-\ga+\frac {d}{2}}\Ga(\tfrac{r}{\ga}+\tfrac {d}{2\ga}).
\end{align*}
 Recalling the formula for $K_r$ from \eqref{Kr}, we obtain
\[ M(z)\sim \frac d2(2\pi)^{-d}\frac{\pi^{\frac d2}}{\Ga(1+\frac d2)}\frac1{\ga-\frac d2-r}\Ga(\tfrac{r}{\ga}+\tfrac {d}{2\ga})z^{\frac d2} = \frac{K_r}{|D|}z^{\frac d2},\qquad z\to\infty. \]
Inserting this into \eqref{Hn} yields the desired result.

If $r=-\frac d2$, we split $H_n(y)$ into $H_{n,1}(y)+H_{n,2}(y)+H_{n,3}(y)$, where
\begin{align*}
H_{n,1}(y)&:=\frac1{\lvert\log \Del\rvert}   \sum_{\la_k\leq \Del^{-1/\ga}} \la_k^{-\frac d2} \phi_k(y)^2,\\
H_{n,2}(y)&:=\frac1{\Del\lvert\log \Del\rvert}  \sum_{\la_k\leq \Del^{-1/\ga}} \la_k^{-\frac d2-\ga} (1-\ee^{-\la_k^\ga\Del}-\la_k^\ga\Del)\phi_k(y)^2,\\
H_{n,3}(y)&:=\frac1{\Del\lvert\log \Del\rvert}  \sum_{\la_k> \Del^{-1/\ga}} \la_k^{-\frac d2-\ga} (1-\ee^{-\la_k^\ga\Del})\phi_k(y)^2.
\end{align*}
Because $|1-\ee^{-x}-x|\leq \frac{x^2}{2}$ for all $x\geq0$, we have from \eqref{est2},
\begin{align*}
H_{n,2}(y)&\lec \frac{\Del}{\lvert\log \Del\rvert} \sum_{\la_k\leq \Del^{-1/\ga}} \la_k^{\ga-\frac d2} \phi_k(y)^2 \lec \frac{1}{\lvert\log \Del\rvert} \to 0
\end{align*}
as $n\to\infty$. Similarly, using \eqref{est1} and the simple bound $1-\ee^{-x}\leq 1$, we obtain
\begin{align*}
H_{n,3}(y)&\lec \frac{1}{\Del\lvert\log \Del\rvert} \sum_{\la_k> \Del^{-1/\ga}} \la_k^{-\frac d2-\ga} \phi_k(y)^2 \lec \frac{1}{\lvert\log \Del\rvert} \to 0.
\end{align*}
Finally, \eqref{est3} and \eqref{Kr} show $H_{n,1}(y)\to (\ga(4\pi)^{\frac d2}\Ga(\frac d2))^{-1}=\frac{K_r}{D}$, which completes the proof of (i).

For (ii), let $\wt V(t,x)$ be defined in \eqref{eA.8}, so that
\[
   H_n(y) = \frac{1}{\tau_n(r)^2} \int_{\la_1-}^\infty t^{r-\ga} (1- \ee^{-t^\ga \Del})) \, \wt V(\dd t,y).
\]
By Abel summation (see \cite[Theorem~4.2]{Apostol76}), we have
\begin{align*}
   H_n(y) &= \frac{1}{\tau_n(r)^2} \lim_{t \to \infty} \Bigl[\wt V(t,y) t^{r-\ga} (1- \ee^{-t^\ga \Del})\Bigr] - \frac{1}{\tau_n(r)^2} \int_{\lambda_1}^\infty \wt V(t,y) \Bigl( t^{r-\ga} (1- \ee^{-t^\ga \Del})\Bigr)'\, \dd t \\
   &\lec \frac{1}{\tau_n(r)^2} \int_{\lambda_1}^\infty (f_{n,1}(t) + f_{n,2}(t))\, \dd t,
\end{align*}
where
\[
   f_{n,1}(t):= (\gamma - r)\, t^{\frac{d}{2}+r-\ga -1} (1- \ee^{-t^\ga \Del}) ,\qquad f_{n,2}(t):= -\ga\, t^{\frac{d}{2}+r-1}\, \ee^{-t^\ga \Del},
\]
and we have used \eqref{eA.8}. The right-hand side no longer depends on $y$, so the integral over $D_0$ contributes a factor $|D_0|$. Using the change of variables $u = t^\ga \Del$, we see that for $i=1,2$,
\[
   \lim_{n \to \infty } \frac{1}{\tau_n(r)^2} \int_{\lambda_1}^\infty f_{n,i}(t)\, \dd t
\]
exists and is finite. This completes the proof of (ii).
\epr

\blem\label{lem:innerprod}
\benu
\item Fix $d \geq 2$. Let $D\subseteq\R^d$ be a bounded connected open set with a piecewise smooth boundary (in the sense of \cite[Definition~1.17]{Hezari18}) and let $D_0$ be an open set with a smooth boundary satisfying $\ov D_0 \subseteq D$. Then there is $C\in(0,\infty)$, only depending on $D$ and $D_0$, such that for all $k> \ell$,
\[ \Biggl\lvert\int_{D_0} \phi_k(x)\phi_\ell(x) \,\dd x\Biggr\rvert \leq C\frac{\la_k^{\frac12}\la_\ell^{\frac18}}{\la_k-\la_\ell}.\]
\item Let $d=1$ and let $D \subseteq\R$ be an interval and let $D_0$ be a subinterval of $D$. Then there is $C\in(0,\infty)$, only depending on $D$ and $D_0$, such that for all $k> \ell$,
\[
  \Biggl\lvert\int_{D_0} \phi_k(x)\phi_\ell(x) \,\dd x\Biggr\rvert \leq C\frac{\la_k^{\frac12}}{\la_k-\la_\ell}.
\]
\eenu
\elem
\bpr For (i),   we can use the eigenfunction properties of $\phi_k$ and $\phi_\ell$ and  Green's second identity to get
\begin{align*} (\la_k-\la_\ell)\int_{D_0} \phi_k(x)\phi_\ell(x)\,\dd x &= \int_{D_0} (\phi_k(x)\Delta\phi_\ell(x)-\phi_\ell(x) \Delta \phi_k(x))\,\dd x\\
&=\int_{\partial D_0} (\phi_k(x)\partial_\nu\phi_\ell(x)-\phi_\ell(x) \partial_\nu \phi_k(x))\,S(\dd x),
\end{align*}
where $S$ is the surface measure and $\partial_\nu$ is the derivative in direction of the outward pointing unit normal field of the boundary $\partial D_0$. Applying the Cauchy--Schwarz inequality, we obtain
\begin{align*} (\la_k-\la_\ell) \Biggl\lvert\int_{D_0} \phi_k(x)\phi_\ell(x)\,\dd x  \Biggl\rvert&\leq \Biggl(\int_{\partial D_0} \phi_k(x)^2\,S(\dd x)\Biggr)^{\frac12}\Biggr(\int_{\partial D_0} (\partial_\nu\phi_\ell(x))^2\,S(\dd x)\Biggr)^{\frac12}\\
&\qquad+\Biggl(\int_{\partial D_0} \phi_\ell(x)^2\,S(\dd x)\Biggr)^{\frac12}\Biggr(\int_{\partial D_0} (\partial_\nu\phi_k(x))^2\,S(\dd x)\Biggr)^{\frac12}.
\end{align*}
Since $d \geq 2$, by \cite[Corollary~1.8]{Hezari18} and \cite[Theorem~1.1]{Christianson15}, it follows that
\begin{align*} (\la_k-\la_\ell) \Biggl\lvert\int_{D_0} \phi_k(x)\phi_\ell(x)\,\dd x  \Biggl\rvert&\lec\la_k^{\frac18}\la_\ell^{\frac12}+\la_\ell^{\frac18}\la_k^{\frac12}\lec \la_\ell^{\frac18}\la_k^{\frac12},
\end{align*}
which implies the assertion of statement (i).

In the case of (ii), since $d=1$, we may as well assume that $D = (0,\pi)$. In this case, $\lambda_k = k^2$ and $\phi_k(x) = \sin(kx)$. Assuming that $D_0 = (y_1,y_2)$, we consider
\[
   \int_{y_1}^{y_2} \sin(k x) \sin(\ell x)\, \dd x.
\]
By additivity of $v \mapsto \int_0^v \cdots$, it suffices to consider $y_1=w_1=0$. The estimate
\[\Biggl\lvert\int_0^{y_2} \sin(ky)\sin(\ell y) \,\dd y\Biggr\rvert=\Biggl\lvert\frac{k\sin(\ell y_2)\cos(ky_2)-\ell\sin(ky_2)\cos(my_2)}{k^2-\ell^2}\Biggr\rvert \leq \frac{2\lambda_k^{\frac12}}{\lambda_k-\lambda_\ell}
\]
completes the proof of the lemma.
\epr

\blem \label{lemapp}
Let $\calg$ be a sub-$\sigma$-field of $\calf$, $I$ be an index set, and for $v\in I$, let $(Y_n(v))_{n \in \N}$ and $Z(v)$ be nonnegative random variables such that
\beq \label{e0app}
   \sup_{v\in I} \sup_{n\in \N} \bbe[Y_n(v)^6] < \infty \qquad\mbox{and}\qquad \sup_{v\in I} \bbe[Z(v)^6] < \infty.
\eeq
Suppose that
\beq \label{e1app}
   \lim_{n \to \infty}\sup_{v\in I} \bbe\biggl[ \Bigl\lvert \bbe[Y_n(v)^2 \mid \calg] - Z(v)\Bigr\rvert\biggr] = 0\qquad{and}\qquad
   \lim_{n \to \infty} \sup_{v\in I} \bbe\Bigl[\var(Y_n(v)^2 \mid \calg)\Bigr] = 0.
\eeq
Then
\beq \label{e21app}
   \lim_{n \to \infty} \sup_{v\in I} \bbe\biggl[\Bigl\lvert Y_n(v)  - \sqrt{Z(v)} \Bigr\rvert\biggr] = 0.
\eeq
\elem

\bpr Recalling that
$$
   \bbe\Bigl[\var(Y_n(v)^2 \mid \calg)\Bigr] = \bbe\biggl[\Bigl(Y_n(v)^2 - \bbe[Y_n(v)^2 \mid \calg]\Bigr)^2\biggr],
$$
we deduce from \eqref{e1app} that
\beq \label{e3app}
   \lim_{n\to \infty} \sup_{v\in I} \Vert Y_n(v)^2 - \bbe[Y_n(v)^2 \mid \calg] \Vert_{L^2(\Om)} = 0.
\eeq
By \eqref{e1app},
$$ 
   \lim_{n\to \infty} \sup_{v\in I} \Vert \bbe[Y_n(v)^2 \mid \calg] - Z(v) \Vert_{L^1(\Om)} = 0.
$$
Because of the bounded moments assumption \eqref{e0app}, we get
\beq \label{e5app}
   \lim_{n\to \infty} \sup_{v\in I} \Vert \bbe[Y_n(v)^2 \mid \calg] - Z(v) \Vert_{L^2(\Om)}= 0.
\eeq
From \eqref{e3app} and \eqref{e5app}, we see that
\beq\label{e6app}
   \lim_{n\to \infty} \sup_{v\in I} \Vert Y_n(v)^2 - Z(v) \Vert_{L^2(\Om)} = 0.
\eeq
Since $Y_n(v)$ and $Z(v)$ are nonnegative,
$$
   \bbe\biggl[\Bigl\lvert Y_n(v) - \sqrt{Z(v)}\Bigr\rvert 1_{\{Y_n(v) + \sqrt{Z(v)} \leq \eps\}}\biggr] \leq \eps.
$$
Therefore,
\begin{align*}
   \bbe\biggl[\Bigl\lvert Y_n(v) - \sqrt{Z(v)}\Bigr\rvert\biggr] &\leq \bbe\biggl[\Bigl\lvert Y_n(v) - \sqrt{Z(v)}\Bigr\rvert 1_{\{Y_n(v) + \sqrt{Z(v)} > \eps\}} \biggr] + \eps\\
   &\leq \bbe\Biggl[\Bigl\lvert Y_n(v) - \sqrt{Z(v)}\Bigr\rvert \frac{Y_n(v) + \sqrt{Z(v)}}{\eps} 1_{\{Y_n(v) + \sqrt{Z(v)} > \eps\}} \Biggr] + \eps \\
   &\leq \frac{1}{\eps} \bbe\biggl[\Bigl\lvert Y_n(v)^2 - Z(v)\Bigr\rvert\biggr]  + \eps.
\end{align*}
From this bound and \eqref{e6app}, we immediately deduce \eqref{e21app}.
\epr

\end{appendix}

\section*{Acknowledgements}
%

The second author was supported in part by the Swiss National Foundation for Scientific Research.


\begin{thebibliography}{4}
	
	\bibitem{Agmon65}
	\textsc{Agmon, S.} (1965). On kernels, eigenvalues, and eigenfunctions of operators related to
	elliptic problems. {\em Comm. Pure Appl. Math.}, \bff{18} 627--663.
	
	\bibitem{Apostol76}
	\textsc{Apostol, T. M.} (1976).
	 {\em Introduction to Analytic Number Theory}.
	 Springer, New York.
	
	\bibitem{BN11} 
	\textsc{Barndorff-Nielsen, O. E., Corcuera, J. M. and Podolskij, M.} (2011).
	 Multipower variation for {B}rownian semistationary processes.
	 \emph{Bernoulli}, \textbf{17} 1159--1194.
	
	
	\bibitem{Basse17}
	\textsc{Basse-O'Connor, A., Lachi{\`e}ze-Rey, R. and Podolskij, M.} (2017).
	 Power variation for a class of stationary increments {L}{\'e}vy
	driven moving averages.
	 {\em Ann. Probab.}, \textbf{45} 4477--4528.
	
	\bibitem{Bibinger19b}
	\textsc{Bibinger, M. and Trabs, M.} (2019).
	 On central limit theorems for power variations of the solution to the
	stochastic heat equation.
	 In {\em
	 	Stochastic Models, Statistics and Their Applications} (A.~Steland, E.~Rafaj{\l}owicz and O.~Okhrin, eds.) 69--84.
	Springer, Cham.
	
	\bibitem{Bibinger19}
	\textsc{Bibinger, M. and Trabs, M.} (2020).
	 Volatility estimation for stochastic {PDE}s using high-frequency
	observations.
	 {\em Stochastic Process. Appl.}, \textbf{130} 3005--3052.
	
	\bibitem{Billingsley95}
	\textsc{Billingsley, P.}  (1995).
	 {\em Probability and Measure},
	 3rd ed., Wiley, New York.
	
	\bibitem{Bingham87}
	\textsc{Bingham, N. H., Goldie, C. M. and Teugels, J. L.}  (1987).
	 {\em Regular Variation}.
	 Cambridge University Press, Cambridge.
	
	\bibitem{Birman80}
	\textsc{Birman, M. \v{S}. and Solomjak, M. Z.} (1980).
	 {\em Quantitative Analysis in {S}obolev Imbedding Theorems and
		Applications to Spectral Theory}.
	 American Mathematical Society, Providence, RI.
	
	\bibitem{Bogachev98}
	\textsc{Bogachev, V. I.} (1998).
	 {\em Gaussian Measures}.
	 American Mathematical Society, Providence.
	
	\bibitem{Chong20}
	\textsc{Chong, C.} (2020).
	 High-frequency analysis of parabolic stochastic {PDE}s.
	 {\em Ann. Statist.}, \textbf{48} 1143--1167.
	
	\bibitem{Christianson15}
	\textsc{Christianson, H., Hassell, A. and Toth, J. A.}  (2015).
	 Exterior mass estimates and ${L}^2$-restriction bounds for {N}eumann
	data along hypersurfaces.
	 {\em Int. Math. Res. Not. IMRN}, \textbf{2015} 1638--1665.
	
	\bibitem{Cialenco20}
	\textsc{Cialenco, I. and Huang, Y.} (2020).
	 A note on parameter estimation for discretely sampled {SPDE}s.
	 {\em Stoch. Dyn.}, \textbf{20} 2050016.
	
	\bibitem{Corcuera06}
	\textsc{Corcuera, J. M., Nualart, D. and Woerner, J. H. C.} (2006).
	 Power variation of some integral fractional processes.
	 \emph{Bernoulli}, \textbf{12} 713--735.
	
	\bibitem{Corcuera09}
	\textsc{Corcuera, J. M., Nualart, D. and Woerner, J. H. C.} (2009).
	 Convergence of certain functionals of integral fractional processes. 
	 {\em J. Theoret. Probab.}, \textbf{22} 856--870. 
	
	\bibitem{Corcuera13}
	\textsc{Corcuera, J. M., Hedevang, E., Pakkanen, M. S. and Podolskij, M.} (2013).
	 Asymptotic theory for {B}rownian semi-stationary processes with
	application to turbulence.
	 {\em Stochastic Process. Appl.}, \textbf{123} 2552--2574.
	
	\bibitem{Dalang98}
	\textsc{Dalang, R. C. and Frangos, N. E.} (1998).
	 The stochastic wave equation in two spatial dimensions.
	 {\em Ann. Probab.}, \textbf{26} 187--212.
	
	\bibitem{Dalang05}
	\textsc{Dalang, R. C. and Sanz-Sol{\'e}, M.} (2005).
	 Regularity of the sample paths of a class of second-order spde's.
	 {\em J. Funct. Anal.}, \textbf{227} 304--337.
	
	\bibitem{Dalang20}
	\textsc{Dalang, R. C. and Sanz-Sol{\'e}, M.} (2020).
	 {\em An Introduction to Stochastic Partial Differential Equations}.
	 In preparation.
	
	\bibitem{DaPrato14}
	\textsc{Da Prato, G. and Zabczyk, J.} (2014).
	 {\em Stochastic Equations in Infinite Dimensions},
	 2nd ed., Cambridge University Press, Cambridge.
	
	\bibitem{Davies74}
	\textsc{Davies, E. B.} (1974).
	 Properties of the {G}reen's functions of some {S}chr{\"o}dinger
	operators.
	 {\em J. Lond. Math. Soc. (2)}, \textbf{7} 483--491.
	
	\bibitem{Elizalde12}
	\textsc{Elizalde, E.} (2012).
	 {\em Ten Physical Applications of Spectral Zeta Functions},
	2nd ed., Springer, Berlin.
	
	\bibitem{Foondun15}
	\textsc{Foondun, M., Khoshnevisan, D. and Mahboubi, P.} (2015).
	 Analysis of the gradient of the solution to a stochastic heat
	equation via fractional {B}rownian motion.
	 {\em Stoch. Partial Differ. Equ. Anal. Comput.}, \textbf{3} 133--158.
	
	\bibitem{Garding53}
	\textsc{G{\aa}rding, L.} (1953).
	 On the asymptotic distribution of the eigenvalues and eigenfunctions
	of elliptic differential operators.
	 {\em Math. Scand.}, \textbf{1} 237--255.
	
	\bibitem{Geisinger14}
	\textsc{Geisinger, L.} (2014).
	 A short proof of {W}eyl's law for fractional differential operators.
	 {\em J. Math. Phys.}, \textbf{55} 011504.
	
	\bibitem{Gilbarg01}
	\textsc{Gilbarg, D. and Trudinger, N. S.} (2001).
	 {\em Elliptic Partial Differential Equations of Second Order},
	 Springer, Berlin.
	
	\bibitem{Giacomin01}
	\textsc{Giacomin, G., Olla, S. and Spohn, H.} (2001).
	 Equilibrium fluctuations for $\nabla{\varphi}$ interface model.
	 {\em Ann. Probab.}, \textbf{29} 1138--1172.
	
	\bibitem{Hawking77}
	\textsc{Hawking, S. W.} (1977).
	 Zeta function regularization of path integrals in curved spacetime.
	 {\em Comm. Math. Phys.}, \textbf{55} 133--148.
	
	\bibitem{Hezari18}
	\textsc{Hezari, H.} (2018).
	 Quantum ergodicity and ${L}^p$ norms of restrictions of
	eigenfunctions.
	 {\em Comm. Math. Phys.}, \textbf{357} 1157--1177.
	
	\bibitem{Hildebrandt19}
	\textsc{Hildebrandt, F. and Trabs, M.} (2019).
	 Parameter estimation for {SPDE}s based on discrete observations in
	time and space.
	 arXiv:1910.01004.
	
	\bibitem{Jacod12}
	\textsc{Jacod, J. and Protter, P.} (2012).
	 {\em Discretization of Processes}.
	 Springer, Berlin.
	
	\bibitem{Jacod03}
	\textsc{Jacod, J. and Shiryaev, A. N.} (2003).
	 {\em Limit Theorems for Stochastic Processes}, 
	2nd ed., Springer, Berlin.
	
	\bibitem{Kallenberg02}
	\textsc{Kallenberg, O.} (2002).
	 {\em Foundations of Modern Probability},
	2nd ed., Springer, New York.
	
	\bibitem{Ledoux91}
	\textsc{Ledoux, M. and Talagrand, M.} (1991).
	 {\em Probability in Banach Spaces}.
	 Springer, Berlin.
	
	\bibitem{McCullagh18}
	\textsc{McCullagh, P.} (2018).
	 {\em Tensor Methods in Statistics},
	 2nd ed., Dover, Mineola.
	
	\bibitem{Nourdin08}
	\textsc{Nourdin,  I.} (2008).
	 Asymptotic behavior of weighted quadratic and cubic variations of
	fractional {B}rownian motion.
	 {\em Ann. Probab.}, \textbf{36} 2159--2175.
	
	\bibitem{Nourdin10}
	\textsc{Nourdin, I., Nualart, D. and Tudor, C. A.} (2010).
	 Central and non-central limit theorems for weighted power variations of fractional {B}rownian motion.
	 {\em Ann. Inst. Henri Poincar\'e Probab. Stat.}, \textbf{46} 1055--1079.
	
	\bibitem{Nualart18}
	\textsc{Nualart, D. and Zeineddine,R.} (2018).
	 Symmetric weighted odd-power variations of fractional {B}rownian motion and applications.
	 {\em Commun. Stoch. Anal.}, \textbf{12} 37–-58.
	
	\bibitem{Peccati11}
	\textsc{Peccati, G. and Taqqu, M. S.} (2011).
	 {\em Wiener Chaos: Moments, Cumulents and Diagrams}.
	 Springer, Milan.
	
	\bibitem{Peszat07}
	\textsc{Peszat, S. and Zabczyk, J.} (2007).
	 {\em Stochastic Partial Differential Equations with L{\'e}vy Noise}.
	 Cambridge University Press, Cambridge.
	
	\bibitem{Pospisil07}
	\textsc{Posp{\'i}{\v s}il, J. and Tribe, R.} (2007).
	 Parameter estimates and exact variations for stochastic heat
	equations driven by space-time white noise.
	 {\em Stoch. Anal. Appl.}, \textbf{25} 593--611.
	
	\bibitem{Prevot07}
	\textsc{Pr{\'e}v{\^o}t, C. and R{\"o}ckner, M.} (2007).
	 {\em A Concise Course on Stochastic Partial Differential Equations}.
	 Springer, Berlin.
	
	\bibitem{SanzSole02}
	\textsc{Sanz-Sol{\'e}, M. and Sarr{\`a}, M.} (2002).
	 H{\"o}lder continuity for the stochastic heat equation with spatially
	correlated noise.
	 In {\em Seminar on Stochastic Analysis, Random Fields and Applications, III} (\textit{Ascona}, 1999) (R. C. Dalang, M. Dozzi and F. Russo, eds.). {\em Progress in Probability} \textbf{52} 259--268. Birkh{\"a}user,
	Basel.
	
	\bibitem{Sheffield07}
	\textsc{Sheffield, S.} (2007).
	 Gaussian free fields for mathematicians.
	 {\em Probab. Theory Related Fields}, \textbf{139} 521--541.
	
	\bibitem{Suquet99}
	\textsc{Suquet, C.} (1999).
	 Tightness in {S}chauder decomposable {B}anach spaces.
	 In {\em Proceedings of the St. Petersburg
		Mathematical Society Volume V}, (N. N.~Uraltseva, ed.) 201--224. American Mathematical
	Society, Providence.
	
	\bibitem{Swanson07}
	\textsc{Swanson, J.} (2007).
	 Variations of the solution to a stochastic heat equation.
	 {\em Ann. Probab.}, \textbf{35} 2122--2159.
	
	\bibitem{vanNeerven08}
	\textsc{van Neerven, J. M. A. M., Veraar, M. C. and Weis, L.} (2008).
	 Stochastic evolution equations in {UMD} {B}anach spaces.
	 {\em J. Funct. Anal.}, \textbf{255} 940--993.
	
	\bibitem{VereJones97}
	\textsc{Vere-Jones, D.} (1997).
	 Alpha-permanents and their applications to multivariate gamma,
	negative binomial and ordinary binomial distributions.
	 {\em New Zealand J. Math.}, \textbf{26} 125--149.
	
	\bibitem{Voros87}
	\textsc{Voros, A.} (1987).
	 Spectral functions, special functions and the {S}elberg zeta
	function.
	 {\em Comm. Math. Phys.}, \textbf{110} 439--465.
	
	\bibitem{Voros92}
	\textsc{Voros, A.} (1992).
	 Spectral zeta functions.
	 In {\em Zeta functions in
		geometry} (N.~Kurokawa and T.~Sunada, eds.) 327--358. Kinokuniya, Tokyo.
	
	\bibitem{Walsh86}
	\textsc{Walsh, J. B.} (1986).
	 An introduction to stochastic partial differential equations.
	 In {\em \'Ecole D'été de Probabilités de Saint-Flour, XIV}---1984. {\em Lecture Notes in Math.} \textbf{1180} 265–439. Springer, Berlin.
	
\end{thebibliography}
\end{document}